\documentclass[12pt]{amsart}
\usepackage{amsmath,amscd}  
\usepackage[pdftex]{hyperref}
\hypersetup{bookmarksopen,linkcolor=blue,citecolor=magenta,colorlinks}
\usepackage[notcite,final]{showkeys}

\numberwithin{equation}{section}

\newtheorem{theorem}[subsection]{Theorem}
\newtheorem{lemma}[subsection]{Lemma}

\newtheorem{corollary}[subsection]{Corollary}

\theoremstyle{definition}

\newtheorem{definition}[subsection]{Definition}

\theoremstyle{remark}

\newtheorem{example}[subsection]{Example}

\newcommand{\F}{{\bf F}}
\newcommand{\field}{{\bf{F}}}

\newcommand{\tr}{\mathop{\rm Tr}}
\newcommand{\Norm}{\mathop{\rm N}}
\newcommand{\LT}{\mathop{\rm LT}}
\newcommand{\LM}{\mathop{\rm LM}}

\newcommand{\binomial}[2]{\genfrac{(}{)}{0pt}{}{#1}{#2}}

\newcommand{\res}{\Delta}
\newcommand{\pol}{\nabla}
\newcommand{\Pol}{{\mathcal P}}
\newcommand{\Res}{{\mathcal R}}
\newcommand{\N}{\mathbb{N}}

\newcommand{\B}{\mathcal{B}}
\newcommand{\PDP}[2]{{\text {PDP}^{#1}_{\leq #2}}}
\newcommand{\IDP}[2]{{\text {IDP}^{#1}_{#2}}}
\newcommand{\pony}{{Lemma~\ref{one trick pony}}}
\newcommand{\fp}{{\bf{F}}_p}
\newcommand{\Pnorm}{\mathop{\rm N}}
\newcommand{\dPnorm}{{\rm N}}
\newcommand{\sltwo}{SL_2(\fp)}

\newcommand{\set}[1]{\{#1\}}
\newcommand{\Id}{\operatorname{1}}

\title[Invariants of $m\,V_2$ Revisited]{Vector invariants for the two dimensional modular representation of a cyclic group of prime order}

\author[Campbell]{H E A Campbell}
\address{Mathematics \& Statistics Department \\
\hfil\break\indent Memorial University of Newfoundland \\ St John's NL A1A 5S7 \\Canada }
\email{eddy@mun.ca}
\author[Shank]{R J Shank}
\address{School of Mathematics, Statistics \& Actuarial Science \\
\hfil\break\indent University of Kent,  Canterbury \\ CT2 7NF, UK}
\email{R.J.Shank@kent.ac.uk}
\author[Wehlau]{D L Wehlau}
\address{Department of Mathematics and Computer Science \\ \hfil\break\indent
       Royal Military College \\ King\-ston, Ontario, Canada \\ K7K 5L0
      }
\email{wehlau@rmc.ca}

\date{\today}


\dedicatory{}

\begin{document}
\begin{abstract}
In this paper, we study the vector invariants of the $2$-dimensional indecomposable representation $V_2$ of the cylic group, $C_p$, of order $p$ over a field $\field$ of characteristic $p$, $\field[m\, V_2]^{C_p}$.  This ring of invariants was first studied by David Richman \cite{richman} who showed that the ring required a generator of degree $m(p-1)$, thus demonstrating that the result of Noether in characteristic 0 (that the ring of invariants of a finite group is always generated in degrees less than or equal to the order of the group) does not extend to the modular case.  He also conjectured that a certain set of invariants was a generating set with a proof in the case $p=2$.  This conjecture was proved by Campbell and Hughes in \cite{campbell-hughes}.  Later, Shank and Wehlau in \cite{cmipg} determined which elements in Richman's generating set were redundant thereby producing a minimal generating set.

We give a new proof of the result of Campbell and Hughes, Shank and Wehlau giving a minimal algebra generating set for the ring of invariants $\field[m\,V_2]^{C_p}$.  In fact, our proof does much more.  We show that our minimal generating set is also a SAGBI basis for $\field[m\,V_2]^{C_p}$.  Further, our results provide a procedure for finding an explicit decomposition of $\field[m\,V_2]$ into a direct sum of indecomposable $C_p$-modules.  Finally, noting that our representation of $C_p$ on $V_2$ is as the $p$-Sylow subgroup of $\sltwo$, we describe a generating set for the ring of invariants $\field[m\, V_2]^{\sltwo}$ and show that $(p+m-2)(p-1)$ is an upper bound for the Noether number, for $m>2$.

\end{abstract}

\maketitle

\tableofcontents

\section{Introduction}\label{intro}
We suppose $G$ is a group represented on a vector space $V$ over a field $\field$.  If $\set{x_1, x_2, \dots, x_n}$ is a basis for the hom-dual, $V^* = \hom_\field(V,\field)$, of $V$, then we denote the symmetric algebra on $V^*$ by
    $$
        \field[V] = \field[x_1, x_2, \dots, x_n]
    $$
and we note that $G$ acts on $f \in \field[V]$ by the rule
    $$
        \sigma(f) (v) = f (\sigma^{-1}(v)).
    $$
As an aside, the notation $\field[V]$ is often used in the literature to denote the ring of regular functions on $V$.  Our notation coincides with the usual notion when the field $\field$ is infinite.  However, for example, if $\field = \field_p$, the prime field, then the functions $x_1$ and $x_1^p$ coincide on $V$.

The ring of functions left invariant by this action of $G$ is denoted $\field[V]^G$.  Invariant theorists often seek to relate algebraic properties of the invariant ring to properties of the representation.  For example, when $G$ is finite of order $|G|$ and the characteristic $p$ of $\field$ does {\em not} divide $|G|$ -- the {\em non-modular} case -- then $\field[V]^G$ is a polynomial algebra if and only if $G$ is generated by reflections (group elements fixing a hyperplane of $V$).  This is a famous result due to Coxeter \cite{Coxeter}, Shephard and Todd \cite{Shephard+Todd}, Chevalley \cite{Chevalley}, and Serre\cite{Serre}.  For another example in the non-modular case, it is known by work of Noether \cite{Noether} (when $p = 0$), Fogarty \cite{Fogarty} and Fleischmann \cite{Fleischmann} (when $p > 0$), that $\field[V]^G$ is generated in degrees less than or equal to $|G|$.  And, in the non-modular case with $G$ finite, it is well-known that $\field[V]^G$ is always Cohen-Macaulay.

The case when $p > 0$, $G$ is finite, $V$ is finite dimensional and $p$ does divide $|G|$ is that of modular invariant theory.  Many results that are well understood in the non-modular case are not yet understood or even within reach in the modular case.  For example, in the modular case it is known that if $\field[V]^G$ is a polynomial algebra then $G$ must be generated by reflections, but this is far from sufficient.  For another example,  in the modular case $\field[V]^G$ is ``most often'' not Cohen-Macaulay.  Finally, in the modular case, there are examples where $\field[V]^G$ requires generators of degrees (much) larger than $|G|$, see below: this paper re-examines the first known such example in considerable detail.

There are now several references for modular invariant theory, see Benson \cite{Benson}, Smith\cite{Smith}, Neusel and Smith\cite{Neusel+Smith}, Derksen and Kemper\cite{Derksen+Kemper}, Campbell and Wehlau\cite{campbell-hughes}.

Invariant theorists also seek to determine generators for $\field[V]^G$ and, if possible, relations among those generators.  A famous example is the case of {\em vector invariants}, see Weyl \cite{Weyl}.  Here we consider the vector space
    $$
        m\, V = \overset{m {\text ~ summands}}{\overbrace{V \oplus V \oplus \dots \oplus V}}
    $$
with $G$ acting diagonally.  The invariants $\field[m\, V]^G$ are called vector invariants, and in this case, we seek to describe, determine or give constructions for, the generators of this ring, a {\em first main theorem for} $(G,V)$.  Once this is done a theorem determining the relations among the generators is referred to as a {\em second main theorem for} $(G,V)$.

The cyclic group $C_p$ has exactly $p$ inequivalent  indecomposable representations over a field $\field$ of characteristic $p$.  There is one indecomposable $V_n$ of dimension $n$ for each $1 \le n \le p$. To see this choose a basis for $V_n$ with respect to which a fixed generator, $\sigma$, of $C_p$ is represented by a matrix in Jordan Normal form.  Since $V_n$ is indecomposable this matrix has a single Jordan block and since $\sigma$ has order $p$ the common eigenvalue must be 1, the only $p^\text{th}$ root of unity in a field of characteristic $p$.  Thus $\sigma$ is represented on $V_n$ by the matrix
  $$
   \begin{pmatrix}
   1 & 0 & 0 & \hdots & 0 & 0\\
   1 & 1 & 0 & \hdots & 0  & 0\\
   0 & 1 & 1 & \hdots & 0  & 0\\
   \vdots & \vdots & \vdots & \ddots & \vdots & \vdots\\
   0 & 0 & 0 & \hdots & 1 & 0\\
   0& 0 & 0 & \hdots  & 1 & 1\\
 \end{pmatrix}.
$$
In order that this matrix have order $p$ (or 1) we must have $n \leq p$.  We call such a basis of $V_n$ for which $\sigma$ is in (lower triangular) Jordan Normal form a {\em triangular} basis.

Observe the following chain of inclusions:
  $$
  V_1 \subset V_2 \subset \cdots \subset V_p.
  $$

If $V$ is any finite dimensional $C_p$-module then $V$ can be decomposed into a direct sum of indecomposable $C_p$-modules:
  $$
  V \cong m_1\, V_1 \oplus m_2\, V_2 \dots \oplus  m_p\, V_p
  $$
where $m_i \in \N$ for all $i$.  This decomposition is far from unique but does have the property that the {\em multiplicities} $m_\ell$ are unique.

We are interested in the representation $m \,V_2$ and the action of $C_p$ on $\field[m\,V_2]$.  The ring of invariants  $\field[m\,V_2]^{C_p}$ was first studied by David Richman \cite{richman}.  He showed that this ring required a generator of degree $m(p-1)$, showing that the result of Noether in characteristic 0 did not extend to the modular case.  He also conjectured that a certain set of invariants was a generating set with a proof in the case $p=2$.  This conjecture was proved by Campbell and Hughes in \cite{campbell-hughes}: the proof is long, complex, and counter-intuitive in some respects.  Later, Shank and Wehlau in \cite{cmipg} determined which elements in Richman's generating set were redundant thereby producing a minimal generating set.

We will show later (and the proof is not difficult), that $\field[m\,V_2]^{C_p}$ is not Cohen-Macaulay, or see Ellingsrud and Skjelbred \cite{vikings}.

In this paper, we give a new proof of the result of Campbell and Hughes, Shank and Wehlau giving a minimal algebra generating set for the ring of invariants $\field[m\,V_2]^{C_p}$.  In fact, our proof does much more.  We show that our minimal generating set is also a SAGBI basis for $\field[m\,V_2]^{C_p}$.  In our view, this result is extraordinary.  Further, our techniques also yield a procedure for finding a decomposition of $\field[m\,V_2]$ into a direct sum of indecomposable $C_p$-modules.

Our paper is organised as follows.  In the second section of our paper, Preliminaries, we provide more details on the the representation theory of $C_p$, our use of graded reverse lexicographical ordering on the monomials in $\field[m\, V_2]^{C_p}$, and define the term {\em SAGBI} basis.  In the next section, Polarisation, we define the polarisation map $\field[V] \to \field[m\, V ]$, its (roughly speaking) inverse, known as restitution, and we note that these maps are $G$-equivariant, hence map $G$-invariants to $G$-invariants.  These techniques allow us to focus our attention on multi-linear invariants.  The next section, Partial Dyck Paths, describes a concept arising in the study of lattices in the plane, see, for example the book by Koshy \cite[p.~151]{Koshy}, and is followed by a section on Lead Monomials.  Here we show that there is a bijection between the set of lead monomials of multi-linear invariants and certain collections of Partial Dyck Paths.  This work requires us to count the number of indecomposable $C_p$ summands in
    $$
        \overset{m}{\otimes} V_2 = \overset{m {\text ~ copies}}{\overbrace{V_2 \otimes V_2 \otimes \cdots \otimes V_2}},
    $$
and in fact we are able to determine a decomposition of $\overset{m}{\otimes} V_2$ as a $C_p$-module, see Theorem~\ref{decomp  theorem}.  Following this, in section \S~\ref{generating set section}, we prove that we have a generating set for our ring of invariants.  The next section describes how our techniques provide a procedure for finding a decomposition of $\field[m\, V_2]$ as a $C_p$-module.  In the final section, noting that our representation of $C_p$ on $V_2$ is as the $p$-Sylow subgroup of $\sltwo$, we are able to describe a generating set for the ring of invariants $\field[m\, V_2]^{\sltwo}$.

We thank the referee for a thorough and careful reading of our paper.

\section{Preliminaries}\label{prelim}
Suppose $\set{e_1, e_2, \dots, e_n}$ is a triangular basis for $V_n$. Note that the $C_p$-module generated by $e_1$ is all of $V_n$. We also note that the indecomposable module $V_n^* = \hom(V_n,\field)$ is isomorphic to $V_n$ since $\dim(V_n^*) = \dim(V_n)$. Because of our interest in invariants we often focus on the $C_p$ action on $V_n^*$ rather than on $V_n$ itself.  Therefore we will choose the dual basis $\set{x_1, x_2, \dots, x_n}$ for $V^*$ to the basis $\set{e_1, e_2,  \dots, e_n}$.   With this choice of basis the matrices representing $G$ are upper-triangular on $V^*$.  We note that $\sigma^{-1}(x_1) = x_1$ and $\sigma^{-1}(x_i )= x_i + x_{i-1}$ for $2 \leq i \leq n$: for convenience, and since $\sigma^{-1}$ also generates $C_p$, we will change notation and write $\sigma$ instead of $\sigma^{-1}$ for the remainder of this paper.  With this convention, we note that $(\sigma-1)^r(x_n) = x_{n-r}$ for $r < n$ and $\dim(V_n)=n$ is the largest value of $r$ such that $x_1 \in (\sigma-1)^{r-1}(V_n^*)$.   We say that the invariant $x_1$ has {\em length} $n$ in this case and write $\ell(x_1)=n$. We observe that the socle of $V_n$ is the line $V_n^{C_p}$ spanned by $\set{e_n}$.  Similarly $(V_n^*)^{C_p}$ has basis $\set{x_1}$.

Note that the kernel of $\sigma-1 : V_i \to V_i$ is $V_i^{C_p}$ which is one dimensional for all $i$.  Thus
 $$
    \dim((\sigma-1)^{j}(V_i))
        = \begin{cases} 0   & \text{if $j-1 \geq i$;}\\
                    i-j & \text{if $j-1 < i$.}
\end{cases}
$$
For $$
  V \cong m_1\, V_1 \oplus m_2\, V_2 \dots \oplus  m_p\, V_p
  $$
this gives $(p-j)m_p + (p-1-j)m_{p-1} + \dots + (i-j)m_i = \dim((\sigma-1)^{j}(V))$ for all $0 \leq j \leq p-1$ and this system of equations uniquely determines the coefficients $m_1, m_2, \dots, m_p$.

Each indecomposable $C_p$-module, $V_n$, satisfies $\dim (V_n)^{C_p} = 1$.  Therefore the number of summands occurring in a decomposition of $V$
  is given by $m_1 + m_2 + \dots + m_p = \dim V^{C_p}$.

Consider $\tr:=\sum_{\tau\in C_p}\tau$, an element of the group ring of $C_p$. If $W$ is any finite dimensional
$C_p$-representation, we also use $\tr$ to denote the corresponding $\field[W]^{C_p}$-module homomorphism,
$$\tr: \field[W]\to \field[W]^{C_p}.$$
Similarly we define
    $$
        \Norm : \F[W] \to \F[W]^{C_p}
    $$
by $\Norm(w) = \prod_{\tau \in C_p} \tau(w)$.  

Note that $(\sigma-\Id)^{p-1} = \sum_{i=0}^{p-1} (-1)^i\binomial{p-1}{i} \sigma^i = \sum_{i=0}^{p-1} \sigma^i = \tr$.  It follows that $\tr(v) = 0$ if $ v \in V_n$ for $n<p$, while $\tr(x_p) = x_1$ in $V_p$.

It is also the case that $V_p \cong \field C_p$ is the only free $C_p$-module and hence also the only projective.

The next theorem plays an important role in our decomposition of $\field[V]_{(d_1,d_2,\dots,d_m)}$ as a $C_p$-module (modulo projectives).  A proof in the case $V= V_n$ may be found in Hughes and Kemper \cite[section 2.3]{hughes-kemper}, and a proof of the version cited here is in Shank and Wehlau \cite[section 2]{nnsub}
\begin{theorem}[Periodicity Theorem]\label{periodicity}
 Let $V = V_{n_1} \oplus V_{n_2} \oplus \dots \oplus V_{n_m}$.  Let $d_1,d_2,\dots,d_m$ be
non-negative integers and write $d_i = q_i p+r_i $ where $0 \leq r_i \leq p-1$ for $i=1,2,\dots,m$.
Then
    $$
    \field[V]_{(d_1,d_2,\dots,d_m)} \cong \field[V]_{(r_1,r_2,\dots,r_m)} \oplus t\,V_p
    $$
as $C_p$-modules for some non-negative integer $t$.
\end{theorem}

Comparing dimensions shows that in the above theorem
    $$
        t = \left( \prod_{i=1}^m \binomial{n_i+d_i-1}{d_i} - \prod_{i=1}^m \binomial{n_i+r_i-1}{r_i}\right) / p\ .
    $$

In this paper, we are primarily interested in the case $V = m\,V_2$.  We denote the basis for the $i^{\text{th}}$-copy of $V_2^*$ in this direct sum by $\set{x_i,y_i}$ and we have $\sigma(x_i) = x_i$ and $\sigma(y_i) = y_i + x_i$.

  For this representation of $C_p$, there is another ``obvious'' family of invariants, namely the
    $$
        u_{ij} = x_iy_j -x_jy_i = \left| \begin{matrix} x_i & y_i \\ x_j & y_j \end{matrix}\right|
    $$
for  $m \ge 2$.

\subsection{Relations involving the $u_{ij}$}\label{uij rels}

We will consider now two important families of relations involving the invariants $u_{ij}=x_iy_j - y_i x_j$. First we consider algebraic dependencies among the $u_{ij}$.   Suppose $m \geq 4$ and let $1\leq i < j < k < \ell\leq m$.  It is easy to verify that $0 = u_{ij} u_{k\ell} - u_{ik}u_{j\ell} + u_{i\ell}u_{jk}$.  It can be shown that these relations generate all the algebraic relations among the $u_{st}$.

It is useful to represent products of the various $u_{st}$ graphically as follows.  We consider the vertices of a regular $m$-gon and label them clockwise by $1,2, \dots,m$.  We represent the factor $u_{ij}$ by an arrow or directed edge from vertex $i$ to vertex $j$.  Thus a product of various $u_{st}$ is represented by a number of directed edges joining the
labelled vertices of the regular $m$-gon.

Returning to the relation $u_{ij} u_{k\ell} - u_{ik}u_{j\ell} + u_{i\ell}u_{jk}$, we say that the middle product in this relation, $u_{ik}u_{j\ell}$, is a {\em crossing} since the arrows representing the two factors $u_{ik}$ and $u_{j\ell}$ cross (intersect).  Rewriting the relation as $u_{ik}u_{j\ell} = u_{ij} u_{k\ell} + u_{i\ell}u_{jk}$, we see that we may replace a crossing with a sum of two non-crossings.  As observed by Kempe \cite{Kempe}, since the length of two (directed) diagonals representing $u_{ik}$ and $u_{j\ell}$ exceeds both the lengths represented by the sides $u_{ij}$ and $u_{k\ell}$ and the two sides $u_{i\ell}$ and $u_{jk}$, we may repeatedly use``uncrossing'' relations to rewrite any product of $u_{st}$'s by a sum of such products without any crossings.   Thus the space generated by the monomials in the $u_{st}$ of degree $d$ has a basis represented by planar directed graphs on $m$ vertices having $d$ directed edges.  Here we allow multiple (or weighted) edges to represent powers such as $u_{ij}^a$ for $a \geq 2$.

Now we consider another important class of relations, this time involving the $u_{st}$ and the $x_r$.  Take $m \ge 3$, let $1\leq i < j < k\leq m$ and consider the matrix
    $$
        \left(   \begin{array}{ccc}
            x_i& x_j & x_k \\
            x_i& x_j & x_k \\
            y_i& y_j & y_k
        \end{array} \right).
    $$
Computing the determinant by expanding along the first row we find $x_i u_{jk} - x_j u_{ik} + x_k u_{ij} = 0$.  Since $x_1, x_2, x_3$ is a partial homogeneous system of parameters in $\field[m\, V_2]$ consisting of invariants it is a partial homogeneous system of parameters in $\field[m\, V_2]^{C_p}$.  The relation $x_1 u_{23} - x_2 u_{13} + x_3 u_{12} = 0$ shows that the product $x_3 u_{12}$ represents 0 in the quotient ring $\field[m\, V_2]^{C_p}/(x_1,x_2)$.  Considering degrees, it is easy to see that $u_{12}$ and $x_3$ do not lie in the ideal of $\field[m\, V_2]^{C_p}$ generated by  $(x_1,x_2)$.  Thus $x_3$ represents a zero divisor in the quotient ring $\field[m\, V_2]^{C_p}/(x_1,x_2)$.  This shows that the partial homogeneous system of parameters $x_1, x_2, x_3$ in $\field[m\, V_2]^{C_p}$ does not form a regular sequence.  Therefore  $\field[m\, V_2]^{C_p}$ is not a Cohen-Macaulay ring when $m \geq 3$.  For $m \leq 2$ the ring of invariants $\field[m\, V_2]^{C_p}$ is Cohen-Macaulay since $\field[V_2]^{C_p}=\field[x_1,\Norm(y_1)]$ is a polynomial ring and $\field[2\, V_2]^{C_p}=\field[x_1,x_2,u_{12},\Norm(y_1),\Norm(y_2)]$ is a hypersurface ring.

Throughout this paper we will use graded reverse lexicographic term orders.  We write $\LM(f)$ for the lead monomial of $f$ and $\LT(f)$ for the lead term of $f$.  We follow the convention that monomials are power products of variables and terms are scalar multiples of power products of variables.  If $S=\oplus_{d=0}^{\infty} S_d$ is a graded subalgebra of a polynomial ring then we say a set $B$ is {\em a SAGBI basis for $S$ in degree $d$} if for every $f \in S_d$ we can write $\LM(f)$ as a product $\LM(f) = \prod_{g \in B} \LM(g)^{e_g}$ where each $e_g$ is a non-negative integer.  If  $B$ is a SAGBI basis for $S$ in degree $d$ for all $d$ then we say that $B$ is a SAGBI basis for $S$.  If $B$ is a SAGBI basis for $S$ then $B$ is an algebra generating set for $S$.  The word SAGBI is an acronym for ``sub-algebra analogue of Gr\"obner bases for ideals", and was introduced by Robbianno and Sweedler in \cite{rs} and (independently) by Kapur and Madlener in \cite{km}.  For a detailed discussion of term orders we direct the reader to Chapter 2 of Cox, Little and O'Shea \cite{closh}.  For a discussion and application of SAGBI bases in modular invariant theory, we recommend Shank's paper \cite{Shank}.

   Given a sequence of variables $z_1,z_2,\dots, z_m$ and an element $E=(e_1,e_2,\dots,e_m)$ we write
   $z^E$ to denote the monomial $z_1^{e_1} z_2^{e_2} \cdots z_m^{e_m}$ and we denote the degree
   $e_1 + e_2 + \dots + e_m$ of this monomial by $|E|$.

   The following well-known lemma is very useful for computing the lead term of a transfer.
   \begin{lemma}  \label{one trick pony}
   Let $t$ be a positive integer. Then
   $$\sum_{i=0}^{p-1} i^t =
              \begin{cases}
                   -1, &\text{if } p-1 \text{ divides } t;\\
                   0,  &\text{if } p-1 \text{ does not divide } t.
              \end{cases}
     $$
   \end{lemma}

For a proof of this lemma see for example, \cite[Lemma~9.4]{coinvars}.

\section{Polarisation}\label{polarisation section}
Let $V$ be a representation of a group $G$ and let $r\in \N$ with $r \geq 2$
  and consider the map of $G$-representations
$$
     \pol^*: r\,V \longrightarrow (r-1)\,V\\
 $$ defined by $\pol^* (v_1,v_2,\dots,v_r) = (v_1 , v_2 , \dots,v_{r-2},v_{r-1} + v_r) $.
We also consider the morphism
$$\res^* : (r-1)\,V \longrightarrow r\,V$$ given by
$\res^*(v_1,v_2,\dots,v_{r-1}) = (v_1,v_2,\dots,v_{r-2},v_{r-1},v_{r-1})$.
Dual to these two maps we have the corresponding algebra homomorphisms:
$$\pol : \F[(r-1)\,V] \longrightarrow \F[r\,V]$$ and
$$\res : \F[r\,V] \longrightarrow \F[(r-1)\,V]$$ given by
$\pol(f) = f \circ \pol^*$ and $\res(F) = F \circ \res^*$.
We also define $\pol_r^* = (\pol^*)^{r-1} : r\,V \to V$ and
$\res_r^* = (\res^*)^{r-1}: V \to r\,V$.

Thus $\pol_r : \F[V] \longrightarrow \F[r\,V]$ is given by
$(\pol_r(f))(v_1,v_2,\dots,v_r) = f(v_1+v_2+\dots+v_r)$
and $\res_r : \F[r\,V] \longrightarrow \F[V]$ is given by
$(\res_r(F))(v) = F(v,v,\dots,v)$.
The homomorphism $\pol_r$ is called {\em (complete) polarisation} and the homomorphism
$\res_r$ is called {\em restitution}.


  The algebra $\F[r\,V]$ is naturally $\N^r$ graded by multi-degree:
  $$\F[r\,V] = \bigoplus_{(\lambda_1,\lambda_2,\dots,\lambda_r) \in \N^r} \F[r\,V]_{(\lambda_1,\lambda_2,\dots,\lambda_r)}$$
  where $$\F[r\,V]_{(\lambda_1,\lambda_2,\dots,\lambda_r)} \cong
     \F[V]_{\lambda_1} \otimes \F[V]_{\lambda_2} \otimes \dots \otimes \F[V]_{\lambda_r}\ .$$
For each multi-degree, $(\lambda_1,\lambda_2,\dots,\lambda_r) \in \N^r$ we have the projection
$\pi_{(\lambda_1,\lambda_2,\dots,\lambda_r)} : \F[r\,V] \to \F[r\,V]_{(\lambda_1,\lambda_2,\dots,\lambda_r)}$.
Given a homogeneous element $f \in  \F[V]$ of total degree $d$, i.e., $f \in \F[V]_d$,
its {\em full polarisation} is the multi-linear function
$\Pol(f) = \pi_{(1,1,\dots,1)}(\pol_d(f)) \in \F[d\,V]_{(1,1,\dots,1)}$.
Thus $\Pol : \F[V]_d \to \F[d\,V]_{(1,1,\dots,1)}$.

 More generally, we may use isomorphisms of the form $\F[V \oplus W] \cong \F[V] \otimes \F[W] $
    to define
    $$\pol_{r_1,r_2,\dots,r_m} = \pol_{r_1}  \otimes  \pol_{r_2} \otimes  \dots  \otimes \pol_{r_m}
         : \F[\oplus_{i=1}^m W_i] \longrightarrow  \F[\oplus_{i=1}^m r_iW_i]\ .$$
  Again, for ease of notation,
   if $f \in \F[\oplus_{i=1}^m W_i]_{(\lambda_1,\lambda_2,\dots,\lambda_m)}$
 we write $\Pol(f) = \pi_{(1,1,\dots,1)}(\pol_{\lambda_1,\lambda_2,\dots,\lambda_m}(f)) \in \F[\oplus_{i=1}^m \lambda_iW_i]_{(1,1,\dots,1)}$.
 Here again we call the multi-linear function $\Pol(f)$ the full polarisation of $f$.

   Similarly we define the restitution map
  $$\res_{r_1,r_2,\dots,r_m} =\res_{r_1}\otimes\res_{r_2} \otimes\dots\otimes\res_{r_m}:
                  \F[\oplus_{i=1}^m r_iW_i] \longrightarrow  \F[\oplus_{i=1}^m W_i]\ .$$
 In this setting, if we have co-ordinate variables such as $x_i,y_i, z_i$ for $W_i$
  we will use the symbols
 $x_{ij},y_{ij},z_{ij}$ with $1 \leq j \leq r_i$ to denote the co-ordinate variables for $r_i W_i$.
 In this notation, restitution is the algebra homomorphism determined by
 $\res_{r_1,r_2,\dots,r_m}(x_{ij}) = x_i$, $\res_{r_1,r_2,\dots,r_m}(y_{ij}) = y_i$,
 $\res_{r_1,r_2,\dots,r_m}(z_{ij}) = z_i$, etc.  For this reason, restitution is sometimes referred
 to as {\em erasing subscripts}.   For ease of notation, we will write $\Res$ to denote
 the algebra homomorphism 
 $\res_{\lambda_1,\lambda_2,\dots,\lambda_m}$
 when restricted to
 $\F[\oplus_{i=1}^m \lambda_i W_i ]_{(1,1,\dots,1)}$.
 Thus if $F \in  \F[\oplus_{i=1}^m \lambda_i W_i]_{(1,1,\dots,1)}$
  then $\Res(F) \in \F[\oplus_{i=1}^m W_i]_{(\lambda_1,\lambda_2,\dots,\lambda_m)}$.
  (However, we will sometimes abuse notation and use $\Res$ to denote $\res_{\lambda_1,\lambda_2,\dots,\lambda_m}$
  when the indices $\lambda_1,\lambda_2,\dots,\lambda_m$ are clear from the context.)

 It is a relatively straightforward exercise to verify that for any $f  \in \F[\oplus_{i=1}^m W_i]_{(\lambda_1,\lambda_2,\dots,\lambda_m)}$ we have $\Res(\Pol(f)) = (\lambda_1!,\lambda_2!,\dots,\lambda_m!) f$.

 Since $\pol^*$ and $\res^*$ are both $G$-equivariant, so are all the homomorphisms $\pol_{r_1,r_2,\dots,r_m}$ and $\res_{r_1,r_2,\dots,r_m}$. In particular, if $f$ is $G$-invariant then so is $\Pol(f)$.  Similarly, $\Res(F)$ is $G$-invariant if $F$ is.  We also note that since the action of $G$ preserves degree an element $f$ is invariant if and only if all its homogeneous components are invariant.

\section{Partial Dyck Paths}\label{Dyck paths section}

In this section we consider a generalization of Dyck paths (see the book by Koshy  \cite[p. 151]{Koshy} for an introduction to Dyck paths). For us, a lattice path is a finite sequence of steps in the first quadrant of the xy-plane starting from the origin.
Each step in the path is given by either the vector (1,0) (an $x$-step) or the vector (0,1) (a $y$-step). The number of
steps in the path is called its {\em length}. The path is said to have height $h$ if $h$ is the largest integer such that the
path touches the line $y=x-h$. A lattice path has {\em finishing height} $h$ if the final step ends at a point on the line $y=x-h$.

Associated to each lattice path of length $d$ is a word of length $d$ , i.e., an ordered sequence of
$d$ symbols each either an $x$ or a $y$.
This word encodes the path as follows:
the $i^\text{th}$ symbol of the word is $x$ if the $i^\text{th}$ step of the path is an $x$-step and the
 $i^\text{th}$ symbol of the word is a $y$ if the $i^\text{th}$ step is a $y$-step.


We will consider two types of lattice paths: (i) partial Dyck paths and (ii) initial Dyck paths of escape height $p-1$.
\begin{definition}
  A {\em partial Dyck path} is a lattice path
  which stays on or below the diagonal (the line with equation $y=x$).  A partial Dyck path of finishing height 0, i.e., which finishes on the diagonal, is called a Dyck path.
\end{definition}

\begin{definition}
  An {\em initial Dyck path of escape height $t$} is a lattice path
  of height at least $t$ and
  which if it crosses above the diagonal does so only after it touches the line $y=x-t$.
  Expressed another way, these are paths which consist of an partial Dyck path of finishing height
  $t$ followed by an entirely arbitrary sequence of $x$-steps and $y$-steps.
\end{definition}

Clearly there are $2^d$ lattice paths of length $d$.  We may associate these paths with the
$2^d$ monomials in $\F[d\,V_2]_{(1,1,\dots,1)} \cong \otimes^d V_2$.
The lattice path $\gamma$ of length $d$ is associated to the word
$\gamma_1 \gamma_2 \cdots \gamma_d$ and is associated to the multi-linear monomial
$\Lambda(\gamma) = z_1 z_2 \cdots z_d$ where
$\begin{cases}
   z_i = x_i, &  \text{if }\gamma_i=x;\\
   z_i = y_i, & \text{if }\gamma_i=y.
\end{cases}$

We let $\PDP{d}{q}$ denote the set of all partial Dyck paths of length $d$ and height at most $q$.
Furthermore, we will denote by $\PDP{d}{q}(h)$ the set of partial Dyck paths of length $d$, height at most $q$ and finishing height $h$.
We write $\IDP{d}{q}$ to denote the set of all initial Dyck paths of escape height $q$ and length $d$.

\section{Lead Monomials}\label{lead monomials section}

We wish to consider the $C_p$-representation $\F[d\,V_2]_{(1,1,\dots,1)}\cong \otimes^d V_2$. We consider a decomposition of  $\otimes^d V_2$ into a direct sum of indecomposable $C_p$-representations. Each summand $V_h$ has a one dimensional socle spanned by an element $f$ and we associate to this summand the monomial $\LM(f)$.  We say that the invariant $f$ has {\em length} $h$ and we write $\ell(f)=h$.  In general a non-zero invariant has length $h$ if $h-1$ is the maximal value of $r$ for which $f$ lies in the image of $(\sigma-\Id)^r$.

In order to study $\F[d\,V_2]_{(1,1,\dots,1)}^{C_p}\cong(\otimes^d V_2)^{C_p}$ we
use the graded reverse lexicographic order determined by $y_1 > x_1 > y_2 > x_2 \dots > y_d > x_d$ and
consider
    $$
        M  = \{\LM(f) \mid f \in (\otimes^d V_2)^{C_p}\}\ .
    $$
We will show that the set map
   $$
    \Lambda :  \PDP{d}{p-2} \sqcup \IDP{d}{p-1} \longrightarrow M
   $$
is a bijection.

   We begin by showing that the image of $\Lambda$ lies inside $M$.  In fact we will show that
   if $\gamma \in \PDP{d}{p-2}(h)$ then $\Lambda(\gamma)$ is the lead monomial of
   an invariant of length $h+1$.  Furthermore if  $\gamma \in \IDP{d}{p-1}$ then $\Lambda(\gamma)$
   is the lead monomial of an invariant of length $p$, i.e, an invariant lying in $\tr(\otimes^d V_2)$.

Consider a path $\gamma \in  \PDP{d}{p-1}(h)$ and let $\gamma_1 \gamma_2 \cdots \gamma_d$ be the
associated word.  We wish to match each symbol $\gamma_j$ which is a $y$ with an earlier symbol
$\gamma_{\rho(j)}$ which is an $x$.  We do this recursively as follows.  Choose the smallest value
$j$ such that $\gamma_j = y$ and for which we have not yet found a matching $x$.  Take $i$ to be
maximal such that $i < j$, $\gamma_i = x$ and $i \neq \rho(s)$ for all $s < j$.  Then we define
$\rho(j) = i$.  Let $I_1 = \{j \mid \gamma_j = y\}$, $I_2 = \rho(I_1)$ and
$I_3 = \{1,2,\dots,d\} \setminus (I_1 \sqcup I_2)$.  Then $|I_1| = |I_2| = (d-h)/2$, $|I_3|=h$ and
$\gamma_i = x$ for all $i \in I_3$.

Define
$$\theta_0(\gamma) = \left(\prod_{j \in I_1} u_{ \rho(j),j}\right) \prod_{i \in I_3} x_i \text{ and }
    \theta_0'(\gamma) = \left(\prod_{j \in I_1} u_{ \rho(j),j}\right) \prod_{i \in I_3} y_i\ .$$
    Then $\theta_0(\gamma) \in (\otimes^d V_2)^{C_p}$ and
   $$
         \LM(\theta_0(\gamma)) =  \prod_{j \in I_1} \LM(u_{ \rho(j),j}) \prod_{i \in I_3} x_i
                  = \prod_{j \in I_1} x_{\rho(j)}y_j \prod_{i \in I_3} x_i = \Lambda(\gamma)\ .
   $$

   \begin{lemma}
     $(\sigma-\Id)^h ( \theta_0'(\gamma)) = h!\,\theta_0(\gamma)$  and thus
    $\ell(\theta_0(\gamma)) \geq h+1$.
   \end{lemma}

   \begin{proof}
     We will prove a more general statement.
      We will show that
      $$
        (\sigma - \Id)^{|E|} (y^E) = |E|!\,  x^E.
      $$
      Note that this also implies that $(\sigma - \Id)^{|E|+1} (y^E) = 0$.
      We proceed by induction on $|E|$.  The result is clear for $|E|=0,1$.
       Assume, without loss of generality, that $e_i \geq 1$ for all $i$ and define $   Z_i \in \N^d$ by $x_i = x^{   Z_i}$.
        For $|E|\geq 2$ we have
     \begin{align*}
         (\sigma -& \Id)^{|E|} (y^E)  = (\sigma - \Id)^{|E|-1}(\sigma-1) (y^E) \\
          & = (\sigma - \Id)^{|E|-1}\left(\sum_i e_i x_i y^{E-   Z_i} +
                      \text{ terms divisible by some }x_k x_\ell \right)\\
          & =  (\sigma - \Id)^{|E|-1}\left(\sum_i e_i x_i y^{E-   Z_i}\right) \\
                      &\quad    \text{since the other terms lie in the kernel of } (\sigma - \Id)^{|E|-1}\\
          & = \sum_i e_i x_i (\sigma - \Id)^{|E|-1}\left(y^{E-   Z_i}\right)\\
          & = \sum_i e_i x_i (|E|-1)!\,  x^{E-   Z_i} \text{ by induction}\\
         & = \sum_i e_i  (|E|-1)!\, x^{E} = \left(\sum_i e_i \right)  (|E|-1)!\, x^E\\
         & = |E|  (|E|-1)!\,  x^E = |E|!\, x^E
     \end{align*}
   \end{proof}

   If $\gamma \in \PDP{d}{p-2}$ then we define $\theta(\gamma) = \theta_0(\gamma)$ and
   $\theta'(\gamma) = \theta'_0(\gamma)$.


   Suppose instead that $\gamma \in \IDP{d}{p-1}$ and let $\gamma_1 \gamma_2 \cdots \gamma_d$
   be the word associated to $\gamma$.  Take $s$ minimal such that the path $\gamma'$ associated to
   $\gamma_1 \gamma_2 \cdots \gamma_s$ is a partial Dyck path of finishing height $p-1$.
   Since $\gamma' \in \PDP{s}{p-1}(p-1)$, from the above we have
     $I_1 = \{ j \leq s \mid \gamma_j = y\}$, $I_2 =  \rho(I_1)$ and
   $I_3 = \{1,2,\dots,s\} \setminus (I_1 \sqcup I_2)$ with
   $|I_1|=|I_2|=(s-p+1)/2$, $|I_3|=p-1$ and $\gamma_i = x$ for all $i \in I_3$.  We further
   define $I_4 = \{i > s \mid \gamma_i = x\}$ and $I_5 = \{i > s \mid \gamma_i = y\}$.
   Define
   $$
   \theta'(\gamma) = \theta_0'(\gamma')  \prod_{i\in I_5} y_i \prod_{i \in I_4} x_i
          =  \prod_{j \in I_1} u_{ \rho(j),j}  \prod_{i \in  I_3 \cup I_5} y_i \prod_{i \in I_4} x_i
   $$
and
   $$
   \theta(\gamma) = 
   \tr\left( \theta_0'(\gamma')\right) \prod_{i\in I_5} y_i \prod_{i \in I_4} x_i
          =  \tr\left(\prod_{i \in  I_3 \cup I_5} y_i\right)\prod_{j \in I_1} u_{ \rho(j),j} \prod_{i \in I_4} x_i
   $$
    Then $\theta(\gamma) \in \tr(\otimes^d V_2) \subset (\otimes^d V_2)^{C_p}$ and
    $\ell(\theta(\gamma))=p$.

       By \pony
   \begin{align*}
     \LM(\theta(\gamma)) & =  \left(\prod_{i \in I_4} x_i \prod_{j \in I_1} \LM(u_{ \rho(j),j})\right) \LM(\tr(\prod_{i\in I_3 \cup I_5} y_i))\\
       & =  \left(\prod_{i \in I_4} x_i \prod_{j \in I_1} x_{ \rho(j)} y_j\right) \prod_{i\in I_3} x_i \prod_{i\in I_5} y_i
        = \Lambda(\gamma)
   \end{align*}

%

     In summary, if $\gamma \in \PDP{d}{p-2}(h)$ then $\theta(\gamma)$ is an invariant
     of length at least $h+1$ and lead monomial $\Lambda(\gamma)$.
     If $\gamma \in \IDP{d}{p-1}$ then $\theta(\gamma)$ is an invariant
     of length $p$ and with lead monomial $\Lambda(\gamma)$.
         Note that since these lead monomials are all distinct, the maps $\theta$ and $\Lambda$ are injective.

     It remains to show that $\Lambda$ is onto $M$ and to determine the exact length of the invariants
     $\theta(\gamma)$ when $\gamma \in \PDP{d}{p-2}$.  We will show that $\Lambda$ is onto by showing
     $|M| = |\PDP{d}{p-2} \sqcup \IDP{d}{p-1}|$.  To determine $|M|$ we examine the number of
     indecomposable summands in the decomposition of $\otimes^d V_2$.

Define non-negative integers $\mu_p^d(h)$ by the direct sum decomposition of the $C_p$-module
$\otimes^d V_2$ over $\F$: $$\bigotimes^d V_2 \cong \bigoplus_{h=1}^p \mu_p^d(h)\, V_h\ .$$

Using the convention $\otimes^0 V_2 = V_1$, we have the following lemma.

\begin{lemma}\label{recursive mu}
Let $p \geq 3$.  Then
  $$
      \mu^0_p(h) = \delta_h^1 \text{ and } \mu^1_p(h) = \delta_h^2,
  $$
and
  $$
  \mu^{d+1}_p(h) =
       \begin{cases}
             \mu^d_p(2),                                 & \text{if } h=1;\\
             \mu^d_p(h-1) + \mu^d_p(h+1), & \text{if } 2 \leq h \leq p-2;\\
             \mu^d_p(p-2),                             & \text{if } h = p-1;\\
             \mu^d_p(p-1) + 2\mu^d_p(p),       & \text{if } h = p;\\
        \end{cases}
  $$
for $d \geq 1$.
\end{lemma}

\begin{proof}
The initial conditions are clear.  The recursive conditions follow immediately from the following three equations which may be found for example in Hughes and Kemper \cite[Lemma~2.2]{hughes-kemper}:
    \begin{align*}
       V_1 \otimes V_2 &\cong V_2\\
       V_h \otimes V_2 &\cong V_{h-1} \oplus V_{h+1} \text{ for all }2\leq h \leq p-1\\
       V_p \otimes V_2 &\cong 2\,V_p.
    \end{align*}
\end{proof}

Next we count lattice paths.
 Let $\nu^d_q(h) = |\PDP{d}{q}(h)|$ for $1 \leq h \leq q$.
  We also define $\bar\nu_q^d = |\IDP{d}{q}|$.   With this notation we have the following lemma.

  \begin{lemma}  \label{recursive nu}
  Let $q \geq 2$.  Then
    $$\nu^0_q(h) = \delta_h^0 \text{ and } \nu^1_q(h) = \delta_h^1\ ,$$
    $$\bar\nu^0_q = 0 \text{ and } \bar\nu^1_q = 0\ ,$$
    and
    $$\nu^{d+1}_q(h) =
        \begin{cases}
               \nu^d_q(1),                                  & \text{if } h=0;\\
               \nu^d_q(h-1) + \nu^d_q(h+1),  & \text{if } 1 \leq h \leq q-1;\\
               \nu^d_q(q-1),                              & \text{if } h = q;\\
        \end{cases}$$  
         and
      $$\bar\nu_q^{d+1} = \nu^d_{q-1}(q-1) + 2\bar\nu^d_q$$ for all $d \geq 1$.
  \end{lemma}

\begin{proof}
All of these equations are easily seen to hold except perhaps the final one.
Its left-hand term $\bar\nu_q^{d+1} = |\IDP{d+1}{q}|$ is the number of initial Dyck paths
of length $d+1$ and escape height $q$.  We divide such paths into two classes: those which first achieve height $q$ on their final step and those which achieve height $q$ sometime during the first $d$ steps.   Paths in the first class are partial Dyck paths of length $d$, height at most $q-1$ and finishing height $q-1$ followed by an $x$-step for the $(d+1)^\text{st}$ step.
  There are $\nu^d_{q-1}(q-1) = |\PDP{d}{q-1}(q-1)|$ such paths.
The second class consists of initial Dyck paths of escape height $q$ and length $d$
followed by a final step which may be either an $x$-step or a $y$-step.  Clearly there are
$2|\IDP{d}{q}|=2\bar\nu^d_q$ paths of this kind.
\end{proof}

  \begin{corollary}\label{counting corollary}
  For all $d \in \N$, all primes $p$
    and all $h=1,2,\dots,p-1$ we have
    $$\mu_p^d(h) = \nu_{p-2}^d(h-1)\quad\text{and}\quad \mu_p^d(p) = \bar\nu^d_{p-1}\ .$$
         \end{corollary}

\begin{proof}
   Comparing the recursive expressions and initial conditions for $\mu_p^d(h)$ and $\nu_{p-2}^d(h-1)$ and for $\mu_p^d(p)$ and $\bar\nu_{p-1}^d$
   given in the previous two lemmas makes the result clear for $p \geq 5$.

   For $p=2$ it is easy to see that
       $\mu_2^d(1)=\nu_0^d(0) = \delta_d^0$ for $d\geq 0$ and
   $\mu_2^d(2) = 2^{d-1} = \bar\nu_1^d$ for $d \geq 1$.

   For $p=3$ and $h=1,2$ we have
   $$
    \mu_3^d(h)  = \nu_1^d(h-1) =
                            \begin{cases}
                                   1,& \text{if } h+d\text{ is odd};\\
                                   0,& \text{if } h+d\text{ is even}.
                             \end{cases}
   $$
   Hence
   $\mu_3^d(3)  = \lfloor \frac{2^d-1}{3} \rfloor$ for $d \geq 0$.
   From the recursive relation
   $\bar\nu_2^{d+1} = \nu^d_{1}(1) + 2\bar\nu^d_2$
    it is easy to see that $\bar\nu_2^d =  \lfloor \frac{2^d-1}{3} \rfloor = \mu_3^d(3)$.
\end{proof}

 This corollary implies that the map $\Lambda$ is a bijection.
  Furthermore for all $d$,
 every element of $\{\LM(f) \mid f \in (\otimes^d V_2)^{C_p}\}$ may be written as a product with
 factors from the set $\{\LM(g) \mid g \in B\}$ where
 \begin{align*}
     B :&= \{x_i \mid 1\leq i \leq d\} \cup \{u_{ij} \mid 1 \leq i < j \leq d\}\\
      &\cup \set{ \tr( \prod_{i=1}^d y_i^{e_i}) \mid 0 \leq e_i \leq 1, \forall\, i=1,2,\dots,d } 
      \ .
  \end{align*}

 We record and extend these results in the following theorem.

 \begin{theorem}\label{decomp  theorem}
   Let $p$ be a prime, let $d \in \N$ and suppose $0\leq h \leq p-2$.  Let $\gamma \in \PDP{d}{p-2} \cup \IDP{d}{p-1}$.
   Then
   \begin{enumerate}
     \item $\LM(\theta(\gamma)) = \Lambda(\gamma)$.
       \label{LM claim}
     \item If $\gamma \in \PDP{d}{p-2}(h)$ then the invariant $\theta(\gamma)$ lies in
              $$\F[d\, V_2]^{C_p}_{(1,1,\dots,1)} \cong (\otimes^d V_2)^{C_p}$$ and has length $h+1$.
        \label{PDP claim}
     \item If $\gamma \in \IDP{d}{p-1}$ then the invariant $\theta(\gamma)$ lies in
               $$\F[d\, V_2]^{C_p}_{(1,1,\dots,1)} \cong (\otimes^d V_2)^{C_p}$$ and has length $p$.
        \label{IDP claim}
     \item  $B$ is a SAGBI basis in multi-degree $(1,1,\dots,1)$ for $\F[d\, V_2]^{C_p}$.
                 \label{FDVSASB}
   \end{enumerate}
   Furthermore, we have the following decomposition of the $C_p$ representation
   $\otimes^d V_2$ into indecomposable summands:
   $$\bigotimes^d V_2 \cong \bigoplus_{\gamma \in \PDP{d}{p-2} \cup \IDP{d}{p-1}} V(\gamma)$$
   where $V(\gamma) \cong V_{h+1}$ is a $C_p$-module  generated by $\theta'(\gamma)$,
   with socle spanned by $\theta(\gamma)$ and
   $$h = \ell(\theta(\gamma)) -1
   = \begin{cases}
    \text{the finishing height of } \gamma;& \text{if } \gamma\in\PDP{d}{p-2}(h);\\
    p-1 & \text{if } \gamma \in \IDP{d}{p-1}.
       \end{cases}$$
 \end{theorem}

\begin{proof}
    The assertions~(\ref{LM claim}) and (\ref{IDP claim}) have already been proved.

  To prove the other assertions we consider the $C_p$-module
  $$ W  = \sum_{\gamma \in \PDP{d}{p-2} \cup \IDP{d}{p-1}} V(\gamma)$$
  generated by the set $\{\theta'(\gamma) \mid \gamma \in \PDP{d}{p-2} \cup \IDP{d}{p-1}\}$.
  The set of vectors $\{\theta(\gamma) \mid \gamma \in \PDP{d}{p-2} \cup \IDP{d}{p-1}\}$
  spanning the socles of the $V(\gamma)$ is linearly independent
  since these vectors have distinct lead monomials.  This implies that the above sum is direct:
  $$ W  = \bigoplus_{\gamma \in \PDP{d}{p-2} \cup \IDP{d}{p-1}} V(\gamma)\ .$$
%
%
   Thus $\dim W =  (\sum_{h=0}^{p-2} \,(h+1)\cdot\nu_p^d(h) ) +  p \cdot \bar\nu^d_p$.
  Applying Corollary~\ref{counting corollary}, yields $\dim W = \dim \otimes^d V_2$.
  Since $W$ is a submodule of  $\otimes^d V_2$ we see that $W = \otimes^d V_2$.
    Furthermore,  any set of (spanning vectors for the) socles in any direct
  sum decomposition of $\otimes^d V_2$ there will be exactly $\nu^d_p(h)$ invariants of
  length $h+1$ for $0 \leq h \leq p-2$ (and $\bar\nu^d_p$ of length $p$).
     Combining this fact with $\ell(\theta(\gamma)) \geq h+1$ for all $\gamma\in \PDP{d}{p-2}(h)$,
  we get $\ell(\theta(\gamma)) = h+1$ for all $\gamma \in \PDP{d}{p-2}(h)$, completing the proof of
  assertion~(\ref{PDP claim}) as well as the final assertion of the theorem.  Assertion(\ref{FDVSASB})
  also follows now since we have
  $\{\LM(f) \mid f \in (\otimes^d V_2)^{C_p}\} = \{\LM(\theta(\gamma)) \mid \gamma \in \PDP{d}{p-2} \cup \IDP{d}{p-1}\}$
  and each of these lead monomials may be factored into a product of lead monomials of elements of $B$.
\end{proof}

\section{A Generating Set}\label{generating set section}
Consider the set
\begin{align*}
\B =  \set{x_i,& \Norm(y_i) \mid 1\leq i \leq m} \cup \set{u_{ij} \mid 1 \leq i < j \leq m}\\
                                 & \cup \set{\tr(y^E) ~|~0 \le e_i \le p-1}\ .
\end{align*}

We will show that $\B$ is a generating set, in fact a SAGBI basis for $\F[m\,V_2]^{C_p}$.
Let $f \in \F[m\,V_2]^{C_p}$ be monic and
multi-homogeneous, of multi-degree$(\lambda_1,\lambda_2,\dots,\lambda_m)$.
Let $A$ denote the subalgebra $\F[\B]$.
We proceed by induction on the total degree $d = \lambda_1+\lambda_2+\dots+\lambda_m$ of $f$. Clearly if $f$ has total degree 0 then $f$ is constant, $f \in A$ and $\LM(f) =1$ and there is nothing more to prove.

 Suppose then that the total degree $d$ of $f$ is positive.  First suppose that $\lambda_i \geq p$
 for some $i$.  We  consider $f$ as a polynomial in $y_i$ and write
 $f = \sum_{j=0}^{\lambda_i} f_j y_i^j$ where $f_j$ is a polynomial which is homogeneous of degree
 $\lambda_i-j$ in $x_i$.  Dividing $f$ by $ \Norm(y_i)$ in $\F[m\,V_2]$ yields $f = q \Norm(y_i) + r$
  where the remainder $r$ is a polynomial whose degree in $y_i$ is at most $p-1$.
  Applying $\sigma$ we have $f = \sigma(f) = \sigma(q)  \Norm(y_i) + \sigma(r)$.
  Since applying $\sigma$ cannot increase the degree in $y_i$, we see that $\sigma(r)$ also has
  degree at most $p-1$ in $y_i$.  By the uniqueness of
  remainders and quotients we must have $\sigma(r)=r$ and $\sigma(q)=q$, i.e.,
  $q,r \in\F[m\,V_2]^{C_p}$.
  Since $\lambda_i \geq p$, we see that $x_i$ divides $r$ and so we have
  $f = q  \Norm(y_i) + x_i r'$ with $q, r' \in \F[m\,V_2]^{C_p}$.
  By induction $q,r' \in A$ and thus $f \in A$.
  Also by induction we have that $\LM(q)$ and $\LM(r')$, hence also $\LM(f)$ may be written
  as products with factors from $\LM(\B)$.

  Therefore, we may assume that $\lambda_i < p$ for all $i=1,2,\dots,m$.
  Then $\kappa = \lambda_1! \lambda_2! \cdots \lambda_m! \neq 0$.
  Define $$F = \Pol(f) \in \F[d\,V_2]_{(1,1,\dots,1)}^{C_p} = (\otimes_{i=1}^d V_2)^{C_p}.$$

   At this point we want to fix some notation.  We will use
   $\{x_{ij},y_{ij} \mid 1 \leq i \leq m, 1 \leq j \leq \lambda_i\}$ as co-ordinate variables for
   $\lambda_1 V_2 \oplus \lambda_2 V_2 \oplus \dots \oplus \lambda_m V_2$.
   We write $u_{ij,k\ell} = x_{ij}y_{k\ell} - x_{k\ell}y_{ij}$.
    We use a graded reverse lexicographic
   order on $\F[\oplus_{i=1}^m \lambda_i V_2]$ after ordering these variables such that the following conditions hold
   \begin{itemize}
        \item $y_{ij} > x_{ij}$,
        \item  if $i<k$ then $y_{ij} > y_{k\ell}$ and $x_{ij} > x_{k\ell}$,
        \item if $j < \ell$ then $y_{ij} > y_{i\ell}$ and $x_{ij} > x_{i\ell}$.
   \end{itemize}

  We will first show that $\B$ generates $\F[m\,V_2]^{C_p}$ as an $\F$-algebra and then show that it is a SAGBI basis.
  Of course, the former statement follows from the latter but we include a separate proof of the former since the proof
  is short and illustrates the main idea we will need for the latter proof.

  By Theorem~\ref{decomp theorem}, we may write
  $$ F = \sum_I \alpha_I \prod_{ij} x_{ij} \prod_{ij,k\ell} u_{ij,k\ell} \prod_E \tr( {\textstyle\prod}_{ij} y_{ij}^{e_{ij}})\ .$$
  Let $e_i = \sum_{j} e_{ij}$.
  \begin{align*}
  f &= \kappa^{-1} \Res(\Pol(f)) = \kappa^{-1} \Res(F)\\
     &=   \kappa^{-1}  \Res\left( \sum_I \alpha_I \prod_{ij} x_{ij} \prod_{ij,k\ell} u_{ij,k\ell} \prod_E \tr( {\textstyle\prod}_{ij} y_{ij}^{e_{ij}})\right)\\
   &=   \kappa^{-1} \sum_I \alpha_I \prod_{ij}\Res(x_{ij}) \prod_{ij,k\ell}\Res(u_{ij,k\ell}) \prod_E\Res(\tr( \prod_{ij} y_{ij}^{e_{ij}}))\\
   &=  \kappa^{-1} \sum_I \alpha_I \prod_{ij} x_i \prod_{ij,k\ell} u_{ik} \prod_E \tr(\prod_i y_i^{e_i}) \in A
  \end{align*}
  where the last equality follows from the following equalities
  \begin{align*}
    \Res(\tr(y^E)) &=\Res(\sum_{\tau \in C_p} \tau(y^E)) = \sum_{\tau \in C_p} \Res(\tau(y^E))
        = \sum_{\tau \in C_p} \tau(\Res(y^E)) \\
        &= \tr(\Res(y^E))\ .
    \end{align*}

  This completes the proof that $\B$ generates $\F[m\,V_2]^{C_p}$ as an $\F$-algebra.
  We continue with the proof that $\B$ is a SAGBI basis.  First we prove a lemma relating our term orders
  and polarisation.

  \begin{lemma}\label{compatible lemma}
    Suppose $\gamma_1, \gamma_2$ are two monomials in $\F[m\,V_2]_{(\lambda_1,\lambda_2,\dots,\lambda_m)}$ with
    $\gamma_1 > \gamma_2$.  Then $\LT(\Pol(\gamma_1)) > \LT(\Pol(\gamma_2))$.
  \end{lemma}
  \begin{proof}
         Write $\gamma_1 = \prod_{i=1}^m x_i^{a_i} y_i^{\lambda_i-a_i}$ and
    $\gamma_2 = \prod_{i=1}^m x_i^{b_i} y_i^{\lambda_i-b_i}$.
          Choose $s$ such that
    $a_s \neq b_s$ but $a_{s+1}=b_{s+1}, \dots, a_m=b_m$.
     Since $\gamma_1 > \gamma_2$ we must have $b_s > a_s$.

    Now
    \begin{align*}
      \LT(\Pol(\gamma_1)) = \prod_{i=1}^m \prod_{j=1}^{a_i} x_{ij} \prod_{j=a_i+1}^{\lambda_i}y_{ij}\text{ and }
      \LT(\Pol(\gamma_2)) = \prod_{i=1}^m \prod_{j=1}^{b_i} x_{ij} \prod_{j=b_i+1}^{\lambda_i}y_{ij}\ .
    \end{align*}
    Writing
    \begin{align*}
       \Gamma_1 &=  \prod_{i=1}^{s-1} \prod_{j=1}^{a_i} x_{ij} \prod_{j=a_i+1}^{\lambda_i}y_{ij},\qquad
       \Gamma_2 =  \prod_{i=1}^{s-1} \prod_{j=1}^{b_i} x_{ij} \prod_{j=b_i+1}^{\lambda_i}y_{ij}\\
      \text{ and} & \qquad \Gamma_0 =  \prod_{i=s+1}^{m} \prod_{j=1}^{a_i} x_{ij} \prod_{j=a_i+1}^{\lambda_i}y_{ij}
    \end{align*}
        we have
  \begin{align*}
        &\LT(\Pol(\gamma_1)) = \Gamma_0\Gamma_1 \prod_{j=1}^{a_s} x_{sj} \prod_{j=a_s+1}^{\lambda_s}y_{sj}\\
        &\text{and} \\
         &\LT(\Pol(\gamma_2)) = \Gamma_0\Gamma_2 \prod_{j=1}^{b_s}x_{sj}\prod_{j=b_s+1}^{\lambda_s}y_{sj}.
 \end{align*}
    Since $a_s < b_s$ we see that $\LT(\Pol(\gamma_1)) > \LT(\Pol(\gamma_2))$.
  \end{proof}

    Write
  $f = \gamma_1 + \gamma_2 + \dots + \gamma_s$ where each $\gamma_i$ is a term
  and $\LM(f) = \LT(f) = \gamma_1$ since $f$ was assumed to be monic.
     Define $F = \Pol(f)$.  By Lemma~\ref{compatible lemma},
  $\LM(F) = \LM(\Pol(\gamma_1))$.
  Furthermore, each monomial of $\Pol(\gamma_1)$ restitutes to $\gamma_1$.
  In particular, $\Res(\Gamma_1) = \gamma_1$ where $\Gamma_1 = \LM(F)$.
  By Proposition~\ref{decomp theorem}(\ref{FDVSASB}),
  we may write
  \begin{align*}
  \Gamma_1 &= \LM(F) =  \LM\left( \prod_{ij} x_{ij} \prod_{ij,k\ell} u_{ij,k\ell} \prod_E \tr( {\textstyle\prod}_{ij} y_{ij}^{e_{ij}})\right)\\
  &= \prod_{ij} x_{ij} \prod_{ij,k\ell} \LM(u_{ij,k\ell}) \prod_E \LM(\tr( \prod_{ij} y_{ij}^{e_{ij}})).
  \end{align*}
  Restituting 
  we find
  \begin{align*}
  \gamma_1&=\Res(\Gamma_1) = \Res\left(\prod_{ij} x_{ij} \prod_{ij,k\ell} \LM(u_{ij,k\ell}) \prod_E \LM(\tr( {\textstyle\prod}_{ij} y_{ij}^{e_{ij}}))\right)\\
  &= \prod_{ij} \Res(x_{ij}) \prod_{ij,k\ell} \Res(\LM(u_{ij,k\ell})) \prod_E \Res(\LM(\tr( {\textstyle\prod}_{ij} y_{ij}^{e_{ij}})))\\
  & =  \prod_{ij} x_{i} \prod_{ij,k\ell} \LM(u_{i,k}) \prod_E \LM(\tr( {\prod}_{i} y_{i}^{\sum_je_{ij}}))
  \end{align*}
  where the last equality follows using Lemma~\ref{LM of traces} below.
  Thus $\LM(f)$ may be written as a product of factors from $\LM(\B)$.  This shows that
  $\B$ is a SAGBI basis for $\F[m\,V_2]^{C_p}$.

  \begin{lemma}\label{LM of traces}
      Let $y^E = \prod_{i=1}^m\prod_{j=1}^{\lambda_i} y_{ij}^{e_{ij}}$ where $e_{ij} \in \{0,1\}$ for all $i,j$.
      Let $e_i = \sum_{j=1}^{\lambda_i} e_{ij}$.
      If $e_i < p$ for all $i=1,2,\dots,m$ then
       $$\Res\left(\LM(\tr( y^E))\right) = \LM\left(\tr(\Res(y^E))\right)\ .$$
  \end{lemma}

  \begin{proof}  Let $s$ be minimal such that $e_1 + e_2 + \dots + e_s \geq p-1$. (If no such $s$
  exists then $\tr(y^E)=0$ and $\tr(\Res(y^E))=0$.)  Let $r$ be minimal such that
  $ e_1 +  e_2 + \dots +  e_{s-1} + e_{s1} + e_{s2} + \dots + e_{sr} = p-1$.
    By \pony
    $$\LM(\tr(y^E)) = \left(\prod_{i=1}^{s-1} \prod_{j=1}^{\lambda_i} x_{ij}^{e_ij}\right)
                                      \prod_{j=1}^r x_{sj}^{e_{sj}} \prod_{j=r+1}^{\lambda_s} y_{sj}^{e_{sj}}
                                       \left(\prod_{i=s+1}^m \prod_{j=1}^{\lambda_i} y_{ij}^{e_{ij}}\right)\ .$$
     Since $\Res(y^E) = \prod_{i=1}^m y_i^{e_i}$, again using
     \pony\ we see that
     $$\LM(\tr(\Res(y^E))) =  \left(\prod_{i=1}^{s-1}  x_i^{e_i}\right)
                                      x_s^t y_s^{e_s-t}
                                       \left(\prod_{i=s+1}^m y_i^{e_i}\right)$$
    where $t = (p-1) - (e_1 + e_2 + \dots + e_{s-1}) = \sum_{j=1}^r e_{ij}$.
    Thus
        $$
            \Res\left(\LM(\tr( y^E))\right) = \LM\left(\tr(\Res(y^E))\right)
        $$
    as required.
  \end{proof}

 \begin{theorem}
 The set
 \begin{align*}
\B' =  \set{x_i,& \Norm(y_i) \mid 1\leq i \leq m} \cup \set{u_{ij} \mid 1 \leq i < j \leq m}\\
                               & \cup \set{\tr(y^E) ~|~0 \le e_i \le p-1,~2(p-1) < |E|}
\end{align*}
is both a minimal algebra generating set and a SAGBI basis for $\F[m\,V_2]^{C_p}$.
 \end{theorem}

 \begin{proof}
 We start by showing $\B'$ is a SAGBI basis.
   We need to see why we do not need invariants of the form
     $\tr( y^E)$ where $|E| \leq 2(p-1)$ as generators.
    To see this, consider such a transfer $\tr(y^E)$.
      By \pony\ its lead term is
      $x_r^{p-1-t+e_r} y_r^{t-p+1} \prod_{i=1}^{r-1} x_i^{e_i}  \prod_{i=r+1}^d y_i^{e_i}$ where
    $r$ is minimal such that $t = \sum_{i=1}^r e_i \geq p-1$.  (We may assume that $r$ exists
    since if $|E| < p-1$ then $\tr( y^E)=0$.)

    Write $\LM(\tr(y^E)) = x_{i_1} x_{i_2} \cdots x_{i_{p-1}} y_{i_p} y_{i_{p+1}} \cdots y_{i_e}$
    where $1 \leq i_1 \leq i_2 \leq \cdots \leq i_e \leq m$.
    Consider $f = \prod_{j=1}^{2p-2-|E|} x_{i_j} \prod_{j=1}^{|E|-(p-1)} u_{i_{p-j},i_{p-1+j}}$.
    Then $\LM(f) = \LM(\tr(y^E))$.   Thus $\{\LM(f) \mid f \in \B'\}$ generates the same algebra as
    $\{\LM(f) \mid f \in \B\}$ which shows that $\B'$ is a SAGBI basis (and hence a generating set)
    for $\F[m\,V_2]^{C_p}$.

    Now we show that $\B'$ is a minimal generating set.  It is clear that the elements $x_i$ and $u_{ij}$
    cannot be written as polynomials in the other elements of $\B'$.  Furthermore, since $\LM(\Norm(y_i))=y_i^p$ is the
    only monomial occuring in any element of $\B'$ which is a pure power of $y_i$, we see that $\Norm(y_i)$
    is required as a generator.  This leaves elements of the form $\tr(y^E)$ with $|E| > 2(p-1)$.
    We proceed similarly to the proof of \cite[Lemma~4.3]{cmipg}.
    Assume by way of contradiction that $\tr(y^E) = \gamma_1 + \gamma_2 + \dots + \gamma_r$ where
    each $\gamma_i$ is a scalar times a product of elements from $\B'\setminus\{\tr(y^E)\}$ and that
    $\LM(\gamma_1) \geq \LM(\gamma_2) \dots \geq \LM(\gamma_r)$.  Then $\LM(\tr(y^E)) \leq \LM(\gamma_1)$.
    First we suppose that $\LM(\gamma_1) = \LT(\tr(y^E))$.  As above we have
    $$
    \LM(\gamma_1) = \LM(\tr(y^E))= x^A y^B
        =x_r^{p-1-t+e_r} y_r^{t-p+1}\prod_{i=1}^{r-1} x_i^{e_i} \prod_{i=r+1}^d y_i^{e_i}
   $$
     where $r$ is minimal such that $t = \sum_{i=1}^r e_i \geq p-1$.

    Since each $e_i < p$ and $\LM(\Norm(y_i)) = y_i^p$ we see that $\Norm(y_i)$ does not divide $\gamma_1$.
But then since $|A| = p-1$ we see that $|A| < |E|-|A|=|B|$ and thus there must be at least one transfer which divides $\gamma_1$.  Conversely since $|A| = p-1$ exactly one transfer (to the first power) may divide $\gamma_1$.  But then the lead monomials of the other factors must divide
$y^B$ and no element of $\B'$ has a lead monomial satisfying this constraint.
This shows that for $|E|>2(p-1)$, the monomial $\LM(\tr(y^E))$ cannot be properly factored using
lead monomials from $\B'$.

Therefore we must have $\LM(\gamma_1) > \LM(\tr(y^E))$ (and $\LM(\gamma_1)=\LM(\gamma_2))$.
Since we may assume that each term of each $\gamma_i$ is homogeneous of degree $E$, we may
write $\LM(\gamma_1) = x^C y^D$ where $C+D=E$.
But $\LM(\tr(y^E))=x^Ay^B$ is the biggest monomial in degree $E$ which satisfies $|A|\geq p-1$.
 Hence $\LM(\gamma_1) > \LM(\tr(y^E))$ implies that $|C| < p-1$.
 Therefore $\gamma_1$ must be a product of elements of the form $x_i, u_{ij}$ and $\Norm(y_i)$ from $\B'$.
As above, since each $e_i < p$, no $\Norm(y_i)$ can divide $\gamma_1$.  But then $\LM(\gamma_1)$ is a product
 of factors of the form $x_i$ and $\LM(u_{ij})=x_iy_j$ and this forces
 $|C| \geq |D| = |E|-|C|$.
 Therefore $2(p-1) > 2|C| \geq |E|$.
 This contradiction shows that we cannot express $\tr(y^E)$ as a polynomial in the other elements of $\B'$
 when $|E| > 2(p-1)$.
\end{proof}

\section{Decomposing $\field[m\, V_2]$ as a $C_p$-module}\label{decomp section}
In this section we show that our techniques give a decomposition of the homogeneous component
    $$
        \field[m\, V_2]_{(d_1,d_2,\dots,d_m)}
    $$
as a $C_p$-module.  We will describe $\field[m\, V_2]_{(d_1,d_2,\dots,d_m)}$ modulo projectives, i.e., we compute the multiplicities of the indecomposable summands $V_k$ of this component for which $k < p$.  Having done this, a simple dimension computation will give the complete decomposition.

By the Periodicity Theorem (Theorem~\ref{periodicity}), we may assume that each $d_i < p$. Let $d=d_1+d_2+\dots+d_m$. The symmetric group on $d$ letters, $\Sigma_d$,  acts on  $\otimes^d V_2$ by permuting the factors.  This action commutes with the action of $C_p$ (in fact with the action of all of $GL(V_2)$). The image of the polarization map consists of those tensors which are fixed by the Young subgroup $Y =\Sigma_{d_1} \times \Sigma_{d_2} \times \dots \times \Sigma_{d_m}$ of $\Sigma_d$.  Since each $d_i < p$, we see that $Y$ is a non-modular group.  Maschke's Theorem then implies that polarization embeds $\field[m\, V_2]_d$ into $\otimes^d V_2$ as a $C_p$-{\em summand}. Therefore $\ell(\Pol(f)) = \ell(f)$ for all $f \in \field[m\, V_2]^{C_p}_{(d_1,d_2,\dots,d_m)}$ and $\ell(\Res(F)) = \ell(F)$ for all $F \in (\otimes^d V_2)^{C_p \times Y}$.

 Using the relations given in Section~\ref{uij rels}, it is straightforward to write down a basis,
consisting of products of $u_{ij}$'s and $x_i$'s, for the invariants in multi-degree $(d_1,d_2,\dots,d_m)$ which lie in the subring
generated by $\{ x_i \mid 1 \leq i \leq m\} \cup \{u_{i,j}\mid 1 \leq i < j \leq m\}$.
Associated to the lead term of each invariant in this basis is an indecomposable summand of $\field[m\, V_2]_{(d_1,d_2,\dots,d_m)}$.
The dimension of this summand may be found using Theorem~\ref{decomp  theorem}. 
More directly, consider a product of $u_{ij}$'s and $x_i$'s, say 
$$f := \prod_{i=1}^m x_i^{a_i} \cdot \prod_{1\leq i < j \leq m} u_{i,j}^{b_{i,j}} \in \field[m\,V_2]^{C_p}\ .$$
It is not too difficult to show that $\LT(f)$ is the lead term of an element of the transfer if and only if there exists $r$ with $1 \leq r \leq m$ such that
$$\sum_{i=1}^r a_i + \sum_{{1 \leq i \leq r \leq j \leq m}\atop {i < j}}^r b_{ij} \geq p-1\ .$$
If no such $r$ exists then $\ell(f) = 1+ \sum_{i=1}^m a_i$ gives the dimension of the associated summand.

Rather than working with the invariants lying in $\field[m\,V_2]$ directly, one may instead use
Theorem~\ref{decomp  theorem} to decompose $\otimes^d V_2$.   It is then possible to perturb this decomposition so that it is a refinement of the splitting given by polarisation/restitution and thus gives a decomposition of 
$\field[m\,V_2]_{(d_1,...,d_m)}$.

\begin{example}
  As an example we compute the decomposition of
    $$
        \field[4\,V_2]_{(p+1,1,1,p+2)}.
    $$
This space has dimension $(p+2)(2)(2)(p+3) = 4p^2 + 20p + 24$.  By Theorem~\ref{periodicity}, we know
    $$
        \field[4\,V_2]_{(p+1,1,1,p+2)} \cong  \field[4\,V_2]_{(1,1,1,2)} \oplus (4p+20)V_p
    $$
and we need to compute the decomposition of
    $$
        \field[4\,V_2]_{(1,1,1,2)}=V_2 \otimes V_2 \otimes V_2 \otimes S^2(V_2).
    $$

We have available the invariants $x_1,x_2,x_3,x_4$ and $u_{12},u_{13},u_{14},u_{23},u_{24},u_{34}$.  Suppose now that $p \geq 7$.  The products of these 10 invariants which lie in degree $(1,1,1,2)$ are as follows (sorted by length):
  \begin{itemize}
  \item[$\ell=2$:]
  $x_4u_{12}u_{34}, \ x_4u_{13}u_{24}, \ x_4u_{14}u_{23}, \ x_1u_{24}u_{34}, \ x_2u_{14}u_{34}, \ x_3u_{14}u_{24}$
  \item[$\ell= 4$:]
    $x_3 x_4^2 u_{12}, \ x_1x_4^2 u_{23}, \ x_1 x_2 x_4 u_{34}, \ x_2 x_4^2 u_{13}, \ x_1x_3 x_4 u_{24}, \ x_2 x_3 x_4 u_{14}$
  \item[$\ell= 6$:] $x_1x_2 x_3 x_4^2$
   \end{itemize}

Consider the invariants of length $2$.  Among the available relations for those of length $2$ we have:
    \begin{align*}
        0 &= x_4(u_{12}u_{34} - u_{13}u_{24} + u_{14}u_{23}), \\
        0 &= u_{34}(x_1 u_{24} -x_2 u_{14} + x_4 u_{23}),\ \text{and}\\
        0 &= u_{14}(x_2 u_{34} -x_3 u_{24} + x_4 u_{23}).
    \end{align*}
Using these three relations we see that the three invariants
    $$
        x_4 u_{13}u_{24}, \quad x_2 u_{14}u_{34}, \quad x_3 u_{14}u_{24}
    $$
may be expressed in terms of the other three invariants
    $$
        x_4 u_{12}u_{34}, \quad x_4 u_{14}u_{23}, \quad x_1 u_{24}u_{34}.
    $$
Furthermore there are no relations involving only these latter three invariants and thus they represent the socles of 3 summands isomorphic to $V_2$.

Among the available relations involving invariants of length 4 we have
   \begin{align*}
        0 &= x_4^2(x_1 u_{23} - x_2 u_{13}  + x_3 u_{12}), \\
        0 &= x_1 x_4(x_2 u_{34} - x_3 u_{24}  + x_4 u_{23}),\ \text{and}\\
        0 &= x_3 x_4(x_1 u_{24} - x_2 u_{14}  + x_4 u_{12}).
   \end{align*}
These allow us to express the three invariants
    $$
        x_2 x_4^2 u_{13}, \quad x_1x_3 x_4 u_{24},  \quad x_2 x_3 x_4 u_{14}
    $$
using only
    $$
        x_3 x_4^2 u_{12}, \quad x_1x_4^2 u_{23}, \quad x_1 x_2 x_4 u_{34}.
    $$
Again these there are no relations involving only these latter 3 invariants and so they represent the socles of 3 summands isomorphic to $V_4$.

Since $x_1x_2 x_3 x_4^2$ spans the socle of a summand isomorphic to $V_6$ we conclude that
    $$
        \field[4\,V_2]_{(1,1,1,2)} \cong 3\,V_2 \oplus 3\,V_4 \oplus V_6 \quad\text{for } p \geq 7.
    $$

For $p=5$, the foregoing is all correct except that the lattice paths corresponding to $x_1 x_2 x_3 x_4^2$ and $x_1 x_2 x_3 x_4 y_4 = \LT(x_1 x_2 u_{34} x_4)$ both attain height $p-1=4$.  Thus in this case these two invariants both represent a projective summand and we have the decomposition
    $$
        \field[4\,V_2]_{(1,1,1,2)} \cong 3\,V_2 \oplus 2\,V_4 \oplus 2\,V_5\quad\text{for } p =5.
    $$

For $p=2,3$ all the relevant lattice paths attain height $p-1$ and so the summand is projective.
Thus
    \begin{align*}
        \field[4\,V_2]_{(1,1,1,2)} &\cong 8\, V_3\quad\text{for } p =3,\text{ and} \\
        \field[4\,V_2]_{(1,1,1,2)} &\cong 12\, V_2\quad\text{for } p =2. \\
    \end{align*}
We will also illustrate how to use the decomposition of $\otimes^5 V_2$ to find the decomposition of $\field[4\,V_2]_{(1,1,1,2)}$ .  By the results of Section~\ref{lead monomials section}, we have $\otimes^5 V_2 \cong 5\,V_2 \oplus 4\,V_4 \oplus V_6$ for $p \geq 7$.
Here the lead monomials are
  \begin{itemize}
  \item[$\ell=2$:]
  $x_1 y_2 x_3 y_4 x_5,\ x_1 x_2 y_3 y_4 x_5,\ x_1 y_2 x_3 x_4 y_5,\ x_1 x_2 y_3 x_4 y_5,\  x_1 x_2 x_3 y_4 y_5$
  \item[$\ell= 4$:]
  $x_1 y_2 x_3 x_4 x_5,\ x_1 x_2 y_3 x_4 x_5,\ x_1 x_2 x_3 y_4 x_5,\ x_1 x_2 x_3 x_4 y_5$
  \item[$\ell= 6$:] $x_1x_2 x_3 x_4 x_5$
   \end{itemize}
and the corresponding invariants are
   \begin{itemize}
  \item[$\ell=2$:]
  $x_5 u_{12}u_{34},\ x_5 u_{14}u_{23},\ x_4 u_{12}u_{35},\ x_1u_{23}u_{45},\ x_1 u_{25} u_{34}$
  \item[$\ell= 4$:]
    $x_3 x_4 x_5 u_{12},\ x_1x_4 x_5 u_{23},\ x_1 x_2 x_5 u_{34},\ x_1 x_2 x_3 u_{45}$
  \item[$\ell= 6$:] $x_1x_2 x_3 x_4 x_5$
   \end{itemize}
   
   The Young subgroup $Y := \Sigma_1 \times \Sigma_1 \times \Sigma_1 \times \Sigma_2$ acts by simultaneously
   interchanging $x_4$ with $x_5$ and $y_4$ with $y_5$.  Clearly the action preserves length.
   The $C_p \times Y$ invariants are
   \begin{itemize}
  \item[$\ell=2$:]
  $x_5 u_{12}u_{34}+x_4 u_{12}u_{35},\ x_5 u_{14}u_{23}+x_4 u_{15}u_{23},\ x_4 u_{12}u_{35}+x_5 u_{12}u_{34},\\
  x_1u_{23}u_{45} + x_1 u_{23} u_{54} = 0,\ x_1 u_{25} u_{34}+ x_1u_{24} u_{35}$
  \item[$\ell= 4$:]
    $x_3 x_4 x_5 u_{12},\ x_1x_4 x_5 u_{23},\ x_1 x_2 x_5 u_{34} + x_1 x_2 x_4 u_{35},\\
     x_1 x_2 x_3 u_{45}+x_1 x_2 x_3 u_{54} = 0$
  \item[$\ell= 6$:] $x_1x_2 x_3 x_4 x_5$
   \end{itemize}

We now restitute these $C_p \times Y$ invariants to $\field[4\,V_2]^{C_p}_{(1,1,1,2)}$. We find
    \begin{align*}
        \Res(x_5 u_{12}u_{34} + x_4 u_{12}u_{35}) &= 2 x_4 u_{12}u_{34},\\
        \Res(x_5 u_{14}u_{23}+x_4 u_{15}u_{23}) & =2x_4 u_{14} u_{23}, \\
        \Res(x_1 u_{25} u_{34} + x_1u_{24} u_{35})& =2x_1u_{24} u_{34}.
    \end{align*}
Thus we find $3$ summands of $\field[4\,V_2]_{(1,1,1,2)}$ isomorphic to  $V_2$.

Restituting the invariants of length 4 we find
    \begin{align*}
        \Res(x_3 x_4 x_5 u_{12}) &= x_3 x_4^2 u_{12},\\
        \Res(x_1x_4 x_5 u_{23})  &= x_1 x_4^2 u_{23},\ \text{and}\\
        \Res(x_1 x_2 x_5 u_{34} + x_1 x_2 x_4 u_{35}) &= 2x_1x_2 x_4 u_{34}.
    \end{align*}
    Thus we have 3 summands isomorphic to $V_4$.
Since $\Res(x_1x_2 x_3 x_4 x_5) = x_1x_2 x_3 x_4^2$, we see that
      $$
        \field[4\,V_2]_{(1,1,1,2)} \cong 3\,V_2 \oplus 3\,V_4 \oplus V_6\quad\text{for } p \geq 7.
      $$
For $p=2,3,5$, the lengths of the above invariants change and we must adjust our conclusions accordingly as we did earlier.  For $p=2$ we must also use the Periodicity Theorem again since $d_4=2=p$.
\end{example}

\section{A First Main Theorem for $\sltwo$}

The purpose of this section is to  use the relative transfer homomorphism to describe the ring of vector invariants, $\field[m\,V_2]^{\sltwo}$. Let $P$ denote the upper triangular Sylow $p$-subgroup of $\sltwo$, giving $\Norm(y) = \Norm^P(y)=y^p-yx^{p-1}$. The ring of invariants of the defining representation of $\sltwo$ is generated by $L = x\Pnorm(y)$ and $D = \Pnorm(y)^{p-1}+x^{p(p-1)}$ (see Dickson \cite{dickson}, Wilkerson \cite{wilkerson}, or Benson \cite[\S8.1]{Benson}). For $\lambda\in{\mathbb N}^m$, define $L_{\lambda} = \pi_{\lambda}\nabla_m(L)$ and $D_{\lambda}=\pi_{\lambda}\nabla_m(D)$, the multi-degree $\lambda$ polarisations. Further define $L_i$ to be the polarisation of $L$ corresponding to $\lambda_i = p+1$ and $\lambda_j=0$ for $j\not=i$. It is easy to verify that $L_i=x_iy^p_i-x_i^py_i$ is the Dickson invariant for the $i^{th}$ summand. Let $L_{ij}$ denote the polarisation corresponding to $\lambda_i=1$, $\lambda_j=p$, and $\lambda_k=0$ otherwise. So, for example, $L_{32}=L_{(0,p,1,0,\ldots,0)}$. Define $${\mathcal D}_m=\left\{\lambda\in{\mathbb N}^m\mid p\ {\rm divides}\ \lambda_i\ {\rm for\ all}\ i\ {\rm and} \sum_{i=1}^m\lambda_i=p(p-1)\right\}. $$ Further define
    \begin{align*}
        {\mathcal S_m} & =\{u_{ij}\mid i<j\leq m\}\cup \{L_i,L_{ij}\mid i,j\in\{1,\ldots, m\},i\not=j\} \\
            &~\quad \cup \{D_{\lambda}\mid\lambda\in{\mathcal D}_m\}.
    \end{align*}

\begin{theorem} \label{sltwogen}
The ring of vector invariants, $\field[ m\,V_2]^{\sltwo}$, is generated
by ${\mathcal S_m}$ and elements from the image of the transfer.
\end{theorem}

Note that the elements of ${\mathcal S}_m$ are clearly $\sltwo$-invariant and
include a system of parameters. Let $A$ denote the
algebra generated by ${\mathcal S} _m$ and let ${\mathfrak a}$ denote the ideal in $\field[ m\,V_2]^P$ generated by
${\mathcal S}_m$.
A basis for the finite dimensional vector space $\field[ m\,V_2]^P/{\mathfrak a}$ lifts to a set of $A$-module generators
for $\field[ m\,V_2]^P$, say ${\mathcal M}$.
Since the relative transfer homomorphism is a surjective $A$-module morphism, $\field[ m\,V_2]^{\sltwo}$ is generated by
${\mathcal S}_m\cup \tr_P^{\sltwo}({\mathcal M})$. The elements of ${\mathcal M}$ may be chosen to be monomials in the
generators of $\field[ m\,V_2]^P$. Since we are working modulo the image of the transfer, it is sufficient to consider monomials
of the form $\Pnorm(y)^{\alpha}x^{\beta}$.

Let ${\mathfrak u}$ denote the ideal  in $\field[ m\,V_2]^P$ generated by $\{u_{ij}\mid i<j\leq m\}$.

\begin{lemma} \label{L-lem}
For $i<j\leq m$, $L_{ij}=x_i\Pnorm(y_j)+u_{ij}x_j^{p-1}$ and $L_{ji}=x_j\Pnorm(y_i)-u_{ij}x_i^{p-1}$, giving
$L_{ij}\equiv_{\mathfrak u}x_i\Pnorm(y_j)$ and $L_{ji}\equiv_{\mathfrak u}x_j\Pnorm(y_i)$.
\end{lemma}
\begin{proof}
Applying $\nabla_m$ to $L$ gives
$$\left(x_1+\cdots + x_m\right)\left(y_1^p+\cdots+y_m^p-(y_1+\cdots +y_m)(x_1+\cdots +x_m)^{p-1}\right).$$
Expanding gives
$$\left(x_1+\cdots + x_m\right)\left(y_1^p+\cdots+y_m^p\right)-(y_1+\cdots +y_m)(x_1+\cdots +x_m)^p.$$
Collecting the appropriate multi-degrees gives $L_{ij}=x_iy_j^p-y_ix_j^p$ and $L_{ji}=x_jy_i^p-y_jx_i^p$.
Using $u_{ij}=x_iy_j-x_jy_i$ and $\Pnorm(y)=y^p-x^{p-1}y$ gives
$$x_i\dPnorm(y_j)+u_{ij}x_j^{p-1}=x_iy_j^p-x_iy_jx_j^{p-1}+x_iy_jx_j^{p-1}-x_j^py_i=L_{ij}$$
and
$$x_j\dPnorm(y_i)-u_{ij}x_i^{p-1}=x_jy_i^p-x_jy_ix_i^{p-1}-y_jx_i^p+x_i^{p-1}x_jy_i=L_{ji}.$$
\end{proof}

Since ${\mathfrak u}\subset {\mathfrak a}$, the preceding lemma and the formula $L_i=x_i\Pnorm(y_i)$ show that
it is sufficient to compute $\tr_P^{\sltwo}$ on monomials of the form $\Pnorm(y)^{\alpha}$ or $x^{\beta}$.

Let $B$ denote the Borel subgroup containing $P$, i.e., the upper triangular elements of $\sltwo$.
Define a weight function on $\field[ m\,V_2]$ by ${\rm wt}(x_i)\equiv_{(p-1)}1$ and ${\rm wt}(y_i)\equiv_{(p-1)}-1$.
Note that $\Pnorm(y_i)$ is isobaric of weight $-1$.
Furthermore, $\field[ m\,V_2]^B$ consists of the span of the the weight zero elements of $\field[ m\,V_2]^P$.
The relative transfer $\tr_P^B$ is determined by weight:
$${\rm Tr}_P^B\left(\Pnorm(y)^{\alpha}x^{\beta}\right)=
\begin{cases}
-\Pnorm(y)^{\alpha}x^{\beta}, & \text{if $(|\beta|-|\alpha|)\equiv_{(p-1)} 0$;}\\
              0 ,         & \text{ otherwise.}
\end{cases}
$$
Since
$\tr_P^{\sltwo}=\tr_B^{\sltwo}\tr_P^B$, it is sufficient to compute $\tr_B^{\sltwo}$ on $\Pnorm(y)^{\alpha}$
with $|\alpha|$ a multiple of $p-1$ and $x^{\beta}$ with $|\beta|$ a multiple of $p-1$. However,
if $|\beta|\geq p-1$, then $x^{\beta}\in \tr^P(\field[ m\,V_2])$ and $\tr_P^{\sltwo}(x^{\beta})\in \tr^{\sltwo}(\field[ m\,V_2])$.
Thus is is sufficient to compute $\tr_B^{\sltwo}(\Pnorm(y)^{\alpha})$ with $|\alpha|$ a multiple of $p-1$.

For $\lambda\in {\mathcal D}_m$, define
$$\dPnorm(y)^{\lambda/p} = \prod_{i=1}^m \dPnorm(y_i)^{\lambda_i/p}.$$

\begin{lemma}\label{D-lem} $\nabla_m(D)\equiv_{\mathfrak u} (\Pnorm(y_1)+\cdots +\Pnorm(y_m))^{p-1}+(x_1+\cdots +x_m)^{p(p-1)}$,
giving $D_{\lambda}\equiv_{\mathfrak u} \binomial{p-1} {\lambda/p}(\Pnorm(y)^{\lambda/p}+x^{\lambda})$ and
$\Pnorm(y)^{\lambda/p}\equiv_{\mathfrak a}-x^{\lambda}$.
\end{lemma}
\begin{proof}
The proof is by induction on $m$.
Note that
$\nabla_m=\nabla^{m-1}$. Thus $\nabla_1(D)=\nabla^0(D)=D=\Pnorm(y)^{p-1}+x^{p(p-1)}$, as required. Recall that the action of
$\nabla$ on $\field[ m\,V_2]$ is determined by $\nabla(x_m)=x_m+x_{m+1}$, $\nabla(y_m)=y_m+y_{m+1}$, $\nabla(x_i)=x_i$,
and $\nabla(y_i)=x_i$ for $i<m$.
Thus $\nabla(u_{ij})=u_{ij}$ if $i<j<m$ and $\nabla(u_{im})=y_i(x_m+x_{m+1})-x_i(y_m+y_{m+1})=u_{im}+u_{i,m+1}$.
Therefore $\nabla$ induces an algebra morphism on $\field[ m\,V_2]^P/{\mathfrak u}$.
Furthermore
$\nabla(\Pnorm(y_i))=\Pnorm(y_i)$ if $i<j<m$ and
$\nabla(\Pnorm(y_m))=y_m^p+y_{m+1}^p-(x_m+x_{m+1})^{p-1}(y_m+y_{m+1})=\Pnorm(y_m)+\Pnorm(y_{m+1})-u_{m,m+1}\sum_{j=0}^{p-2}(-x_m)^jx_{m+1}^{p-2-j}$.
By induction,
    \begin{align*}
        \nabla_{m+1}(D) & =\nabla(\nabla_m(D)) \in \nabla((\dPnorm(y_1)+\cdots +\dPnorm(y_m))^{p-1}\\
            &+(x_1+\cdots +x_m)^{p(p-1)}+{\mathfrak u}).
        \end{align*}
Evaluating the algebra morphism $\nabla$ gives
    \begin{align*}
        \nabla_{m+1}(D) &\in (\nabla(\dPnorm(y_1))+\cdots +\nabla(\dPnorm(y_m)))^{p-1} +(x_1+\cdots +x_{m+1})^{p(p-1)} \\
        & \qquad +\nabla({\mathfrak u})\\
        &\in (\dPnorm(y_1)+\cdots +\dPnorm(y_{m+1}))^{p-1}+(x_1+\cdots +x_{m+1})^{p(p-1)}+{\mathfrak u},
\end{align*}
as required.
\end{proof}

Using the lemma, if $p-1$ divides $|\alpha|$ then $\tr_B^{\sltwo}(\Pnorm(y)^{\alpha})$ is decomposable modulo the image of the transfer, completing the proof of Theorem~\ref{sltwogen}

To complete the calculation of a generating set for $\field[ m\,V_2]^{\sltwo}$ and compute an upper bound for the Noether number, we need only identify a set of $A$-module generators for $\field[ m\,V_2]$. This can be done by applying the Buchberger algorithm to ${\mathcal S}_m$. For example, a Magma  \cite{magma} calculation for $m=3$ and $p=3$, produces $522$ $A$-module generators giving rise to $74$ non-zero elements in the image of the transfer. Subducting the transfers against ${\mathcal S}_m$ gives $11$ new generators and $29$ in total. Magma's {\it MinimalAlgebraGenerators} command reduces the number of generators to $28$, occuring in degrees $2$, $4$, $6$ and $8$. The same calculation for $p=5$ and $m=3$ gives a Noether number of $24$. Thus for $p\in\{3,5\}$ and $m=3$, the Noether number is $(p+m-2)(p-1)=(p+1)(p-1)$.

\begin{theorem} $\field[ m\,V_2]^{\sltwo}$ is generated  as an $A$-module in degrees less than or equal to $(p+m-2)(p-1)$.
\end{theorem}
\begin{proof}
Define ${\mathfrak a}'$ to be the ideal in $\field[ m\,V_2]$ generated by ${\mathcal S}_m$, i.e.,
    $$
        {\mathfrak a}' = A^+\field[ m\,V_2].
    $$
A basis for $\field[ m\,V_2]/{\mathfrak a}'$ lifts to a set of $A$-module generators for $\field[ m\,V_2]$.
We may choose the $A$-module generators to be monomials, $y^{\alpha}x^{\beta}$, which are minimal
representatives of their mod-${\mathfrak a}'$
congruence class. For convenience, denote $d=|\alpha|+|\beta|$.
For $i<j$, using $u_{ij}=x_iy_j-x_jy_i$, if $x_i$ divides  $y^{\alpha}x^{\beta}$,
then $y_j$ does not. For $j\leq i$, using $L_i$ and
$L_{ij}$, if $x_i$ divides  $y^{\alpha}x^{\beta}$, then $y_j^p$ does not. The remaining representatives fall into
two classes: $y^{\alpha}$ and
$y_1^{\alpha_1}\cdots y_k^{\alpha_k}x_k^{\beta_k}\cdots x_m^{\beta_m}$ with $\beta_k\not=0$ and $\alpha_i\leq p-1$.

Case 1: $y^{\alpha}$. Using $D_{\lambda}$ with $\lambda\in{\mathcal D}_m$, we see that, for
$|\gamma|\geq p-1$, $(y^{\gamma})^p$ does not divide $y^{\alpha}$.
Write $\alpha_i=q_ip+r_i$ with $r_i<p$. Then $y^{\alpha}=(y^q)^py^r$ with $|q|\leq p-2$. Thus
$|\alpha|=p|q|+|r|\leq p(p-2)+m(p-1)=(p+m-1)(p-1)-1$. However,
$\tr(y^\alpha)=0$ unless $p-1$ divides $|\alpha|$. Therefore, the $A$-module generators
of the form $\tr(y^{\alpha})$ satisfy $d=|\alpha|\leq (p+m-2)(p-1)$.

Case 2: $y_1^{\alpha_1}\cdots y_k^{\alpha_k}x_k^{\beta_k}\cdots x_m^{\beta_m}$ with $\beta_k\not=0$ and $\alpha_i\leq p-1$.
For $i<j$, let $x^{\gamma}$ be a monomial in $x_1,\ldots,x_{j-1}$. If $|\gamma|=p-1$, then
$x^{\gamma} L_{ij}=x^{\gamma}(x_iy_j^p-x_j^py_i)\equiv_{\mathfrak u} y^{\gamma}y_ix_j^p-x^{\gamma}y_ix_j^p$.
Therefore, if $\beta_j\geq p$ for any $j>k$, then $|\alpha|<p$.
If $|\gamma|\leq p-1$ then
$x^{\gamma} L_{j}=x^{\gamma}(x_jy_j^p-x_j^py_j)\equiv_{\mathfrak u} y^{\gamma}y_j^{p-|\gamma|}x_j^{|\gamma|+1}-x^{\gamma}y_jx_j^p$.
Therefore, if $\beta_k\geq p$, we also have $|\alpha|<p$. If $|\beta_j|<p$ for all $j\geq k$, then
$|\alpha|+|\beta|\leq (m-k+2)(p-1)\leq (p+m-2)(p-1)$ if $k>1$. Hence it is sufficient to consider the case $|\alpha|<p$.
However the transfer is zero unless $p-1$ divides $|\alpha|$ so we may assume $|\alpha|=p-1$.
If $|\alpha|=p-1$, a straightforward calculation with binomial coefficients gives
$\tr^P(y^{\alpha}x^{\beta})=-x^{\alpha+\beta}$. Furthermore, $\tr_P^B(x^{\alpha+\beta})=0$ unless $p-1$ divides $|\alpha|+|\beta|$.
Write $\alpha_i+\beta_i=q_ip+r_i$ with $r_i<p$. Then $x^{\alpha+\beta}=(x^q)^px^r$.
If $|q|\geq p-1$ and $|r|>0$, we may choose $i$ so that $r_i>0$, choose $\lambda\in {\mathcal D}_m$
so that $x^{\lambda}$ divides $x^{pq}$ and choose $j$
so that $x_j$ divides $x^{\lambda}$. By Lemma~\ref{L-lem}, $x_i\Norm(y_j)\in{\mathfrak a}'$.
Form the S-polynomial between $D_{\lambda}$ and $x_i\Norm(y_j)$.
Using Lemma~\ref{D-lem}, this S-polynomial reduces to $x_i x^{\lambda}$.
Thus either $|q|<p-1$ or $|q|=p-1$ and $|r|=0$.
If $|r|=0$ and $|q|=p-1$, then $d=(p-1)(p-1)\leq (p+m-2)(p-1)$.
Suppose $|q|<p-1$. Then $d\leq m(p-1)+p(p-2)= (p+m-1)(p-1)-1$.
Since $d$ must be a multiple of $p-1$, we have $d\leq (p+m-2)(p-1)$.
\end{proof}

\begin{corollary} For $m>2$, the Noether number for $\field[ m\,V_2]^{\sltwo}$ is less than or equal to
$(p+m-2)(p-1)$. For $m=2$ and $p>2$, the Noether number is $p(p-1)$ and for $m=2$, $p=2$, the Noether
number is $p+1=3$.
\end{corollary}
\begin{proof} The elements of ${\mathcal S}_m$  lie in degrees $2$, $p+1$ and $p(p-1)$.
Clearly $L_1$ and $D_{(p(p-1),0,\ldots,0)}$ are indecomposable.
\end{proof}

For $p=2$ and  $m\in\{3,4\}$, Magma calculations give the Noether number $(p+m-2)(p-1)=m$.


\ifx\undefined\bysame
\newcommand{\bysame}{\leavevmode\hbox to3em{\hrulefill}\,}
\fi


\begin{thebibliography}{30}

\bibitem{Benson}
D.~J.~Benson, {\em Polynomial Invariants of Finite Groups}, Lond.\ Math.\ Soc.\ Lecture Note Ser. {\bf 190} (1993),
Cambridge Univ.\ Press.

\bibitem{magma}
W.~Bosma, J.~J.~Cannon and C.~Playoust,
{\em The Magma algebra system I: the user language},
J.\ Symbolic Comput.\ {\bf 24} (1997), 235--265.

\bibitem{campbell-hughes}
H.~E.~A.~Campbell and I.~P.~Hughes, {\em Vector invariants of $U_2({\bf F}_p)$: A proof of a
conjecture of Richman}, Adv.\ in Math.\ {\bf 126} (1997), 1--20.

\bibitem{campbell-wehlau}
H.~E.~A.~Campbell and David~L.~Wehlau, {\em Modular Invariant Theory}, to appear.

\bibitem{coinvars}
H.~E.~A.~Campbell, I.~P.~Hughes, R.~J.~Shank and D.~L.~Wehlau, {\em Bases for rings of coinvariants},
Transformation Groups {\bf 1} (4) (1996), 307--336.

\bibitem{Chevalley}
C.~Chevalley, {\em Invariants of finite groups generated by reflections},
  Amer.\ J.\ Math.\ \textbf{77} (1955), 778--782.

\bibitem{closh}
D.~Cox, J.~Little, and D.~O'Shea, {\em Ideals, varieties, and algorithms}, (1992)
Springer-Verlag.

\bibitem{Coxeter}
H.~S.~M. Coxeter, \emph{The product of the generators of a finite group
  generated by reflections}, Duke Math. J. \textbf{18} (1951), 765--782.

\bibitem{Derksen+Kemper}
H.~Derksen and G.~Kemper, {\em Computational invariant theory}, Invariant Theory and Algebraic Transformation Groups, I, {\bf 130} (2002), Encyclopaedia of Mathematical Sciences, Springer-Verlag.

\bibitem{dickson}
L.~E.~Dickson,
{\em A Fundamental System of Invariants of the General Modular Linear Group with
a Solution of the Form Problem},
Trans.\ Amer.\ Math.\ Soc. {\bf 12} (1911), 75--98.

\bibitem{vikings}
G.~Ellingsrud and T.~Skjelbred, {\em Profondeur d'anneaux d'invariants en caract{\'e}ristique {$p$}}, Compositio Math.\ {\bf 41} No. 2 (1980), 233--244.

\bibitem{Fogarty}
J.~Fogarty,  {\em On Noether's bound for polynomial invariants of a
finite group},  Electron.\ Res.\ Announc.\ Amer.\ Math.\ Soc.\ { 7} (2001), 5--7.

\bibitem{Fleischmann}
P.~Fleischmann, {\em The Noether bound in invariant theory of
finite groups}, Adv.\ in Math.\ {\bf 152} (2000) no.~1, 23--32.

\bibitem{hughes-kemper}
I.~P.~Hughes and G.~Kemper, {\em Symmetric powers of modular representations,
Hilbert series and degree bounds}, Comm.\ in Alg.\ {\bf 28} (2000), 2059--2088.

\bibitem{km}
D.~Kapur and K.~Madlener, {\em A completion procedure for computing a canonical basis of a
k-subalgebra}, Proceedings of Computers and Mathematics {\bf 89} (1989), ed. E. Kaltofen and S. Watt, MIT, 1--11.

\bibitem{Kempe}
A.~Kempe, {\em On regular difference terms}, Proc. London Math. Soc. 25 (1894), 343-350.

\bibitem{Koshy} T.~Koshy, {\em Catalan Numbers with Application} Oxford University Press, (November 2008).

\bibitem{Neusel+Smith}  M.~D.~Neusel and L.~Smith, {\em Invariant theory of finite groups},
 American Mathematical Society, Providence, RI, 2002, Mathematical Surveys and Monographs, {\bf 94}.

\bibitem{Noether} E.~Noether,  {\em Der Endlichkeitssatz der invarianten endlicher
Gruppen}, Math.\ Ann.\ {\bf 77}, 1915, 89--92; reprinted in: Collected
Papers, Springer-Verlag, Berlin, 1983, 181--184.

\bibitem{rs}
L.~Robbianno and M.~Sweedler, {\em Subalgebra bases}, Lecture Notes in Math.\ {\bf 1430} Springer (1990), 61--87.

\bibitem{richman} D.~R.~Richman,
{\em On vector invariants over finite fields},
Adv.\ in Math.\ {\bf 81} (1990) no.1, 30--65.

\bibitem{Serre}
J.-P.~Serre, {\em Groupes finis d'automorphismes d'anneaux locaux
  r{\'e}guliers}, Colloque d'Alg{\`e}bre (Paris, 1967), Exp. 8, Secr{\'e}tariat
  math{\'e}matique, Paris, 1968, p.~11.

\bibitem{Shank} R.~J.~Shank,
{\em S.A.G.B.I. bases for rings of formal modular seminvariants},
Comment.\ Math.\ Helv.\ {\bf 73} (1998) no.~4, 548--565.

\bibitem{cmipg} R.~J.~Shank and D.~L.~Wehlau,
  {\em Computing modular invariants of {$p$}-groups},
  J.\ Symbolic Comput.\ {\bf 34} (2002) no. 5, 307--327.

\bibitem{nnsub} R.~J.~Shank and D.~L.~Wehlau,
  {\em Noether numbers for subrepresentations of cyclic groups of prime order},
  Bull. London Math. Soc. \ {\bf 34} (2002) no. 4, 438--450.

\bibitem{Shephard+Todd}G.~C.~Shephard and J.~A.~Todd, {\em Finite unitary reflection groups},
  Canadian J.\ Math.\ {\bf 6} (1954), 274--304.

\bibitem{Smith}
L.~Smith, {\em Polynomial invariants of finite groups}, Research Notes in
  Mathematics, vol.~6, A K Peters Ltd., Wellesley, MA, 1995.

\bibitem{Weyl}
H.~Weyl, {\em The classical groups}, Princeton University Press, 1997.

\bibitem{wilkerson} C.~W.~Wilkerson,
{\em A Primer on the {D}ickson Invariants},
Amer.\ Math.\ Soc.\ Contemp.\ Math.\ Series {\bf 19} (1983), 421--434.
\end{thebibliography}
\end{document}